\documentclass{article}

\usepackage{amsmath,amssymb,epsf,wrapfig}
\usepackage[dvipdfm]{graphicx}
\def\qed{\hfill$\square$}
\newtheorem{theorem}{Theorem}
\newtheorem{corollary}[theorem]{Corollary}
\newtheorem{lemma}[theorem]{Lemma}
\newtheorem{proposition}[theorem]{Proposition}

\newtheorem{remark}[theorem]{Remark}

\newtheorem{definition}[theorem]{Definition}

\begin{document}

\title{Chart description for hyperelliptic Lefschetz fibrations and their stabilization}
\author{Hisaaki Endo\\ 
Department of Mathematics, Tokyo Institute of Technology, \\
2-12-1 Oh-okamaya, Meguro-ku, Tokyo 152-8551, Japan \\
\\
Seiichi Kamada\\
Department of Mathematics, Osaka City University, \\
3-3-138 Sugimoto, Sumiyoshi-ku, Osaka 558-8585, Japan}

\date{}

\maketitle

\begin{abstract}
Chart descriptions are a 
graphic method to describe monodromy representations of various topological objects. Here we introduce a chart description for hyperelliptic Lefschetz fibrations, and show that any hyperelliptic Lefschetz fibration can be stabilized by fiber-sum with certain basic Lefschetz fibrations. 
\end{abstract}

\section{Introduction}

Chart descriptions were originally introduced in order to describe $2$-dimensional braids  in \cite{Kam92, Kam96} (cf. \cite{Kam02}).  
In \cite{KMMW}, a chart description for genus-one Lefschetz fibrations was introduced and an elementary proof of Matsumoto's  classification theorem was given.  
At the third JAMEX meeting in Oaxaca, Mexico, 2004, the second author generalized it to a method 
describing any monodromy representation \cite{Kam07} 
and investigated genus-two Lefschetz fibrations as an application \cite{Kam12}. 
Here we introduce a chart description for hyperelliptic Lefschetz fibrations, and show that any hyperelliptic Lefschetz fibration can be stabilized by fiber-sum with certain basic Lefschetz fibrations.

\section{Lefschetz fibrations} 

Let $M$ and $B$ be compact, connected, and 
oriented smooth $4$-manifold and $2$-manifold, respectively.   
Let $f : M \to B$ be a smooth map with $\partial M = f^{-1}(\partial B)$.    
A critical  point $p$ is called a {\it Lefschetz singular point} of 
{\it positive type} (or of {\it negative type}, respectively) 
if there exist local complex 
coordinates $z_1, z_2$ around $p$ and 
a local complex coordinate $\xi$ around $f(p)$  
such that $f$ is locally written as $\xi= f(z_1,z_2) =z_1z_2$ (or $\overline{z_1}z_2$, resp.).  
We call  $f$ 
a  (smooth or differentiable) {\it  Lefschetz fibration\/} if 
all critical points are Lefschetz singular points and if 
there exists exactly one critical point in the preimage of each critical value.  

A {\it general fiber\/} is the preimage of a regular value of $f$.  
The {\it genus} of a Lefschetz fibration is defined to be the genus $g$ of a general fiber. 
A {\it singular fiber\/} of {\it positive type} (or  {\it negative type}, resp.) is the preimage 
of a critical value which contains a Lefschetz singular point of positive type (or negative type, resp.).  
A singular fiber is obtained by shrinking a simple loop, called a vanishing cycle, on a general fiber.  
In this paper we assume that a Lefschetz fibration is `{relatively minimal}', i.e., all vanishing cycles are essential loops.  
We say that a singular fiber is of  {\it type {\rm I}} if the vanishing cycle is a non-separating loop. We say that a singular fiber is 
of {\it type ${\rm II}_h$} for $h=1,\ldots , [g/2]$ if the vanishing cycle is 
a separating loop which bounds a genus-$h$ subsurface of the general fiber.  

A singular fiber is of {\it type ${\rm I}^+$} if it is of type I and of positive type.  Similarly 
{\it type ${\rm I}^-$} and {\it type ${\rm II}_h^+$},  {\it type ${\rm II}_h^-$} 
for $h=1,\ldots ,[g/2]$ are defined.   
We denote by 
$n_0^{+}(f)$, $n_0^{-}(f)$, $n_h^{+}(f)$, and $n_h^{-}(f)$, 
the numbers of singular fibers of $f$ of type ${\rm I}^+$, ${\rm I}^-$, ${\rm II}_h^+$, and ${\rm II}_h^-$, respectively.  
A Lefschetz fibration is called  {\it irreducible} if 
every singular fiber is of type I, i.e.,  $n_h^+(f) = n_h^-(f) = 0$ for $h=1,\ldots ,[g/2]$.   
A Lefschetz fibration is called {\it chiral} or {\it symplectic} if every singular fiber is of positive type, i.e.,  $n_0^-(f)=n_h^-(f) =0$ for $h=1,\ldots ,[g/2]$.

Let $f : M \to B$ be a Lefschetz fibration, and  
$\Delta =\{q_1, \dots, q_n\}$ the set of critical values.  
Let $\rho: \pi_1(B \setminus \Delta, q_0) \to MC$ be the monodromy representation of $f$, where $q_0$ is a base point of $B \setminus \Delta$ and $MC$ is the mapping class group of the fiber $f^{-1}(q_0)$.  
Consider a Hurwitz arc system for 
$\Delta$, say ${\cal A}= (A_1, \dots, A_n)$; each $A_i$ is an embedded arc in $B$ connecting $q_0$ 
and a point of $\Delta$ such that $A_i \cap A_j= \{q_0\}$ for $i \ne j$, and they appear in this order around $q_0$.  
When $B$ is a $2$-sphere or a $2$-disk, the system ${\cal A}$ determines a system of generators of $\pi_1(B \setminus \Delta, q_0)$, say $(a_1, \dots, a_n)$.  
We call $( \rho(a_1), \dots, \rho(a_n) )$ a {\it Hurwitz system} of $f$.  
For details on Hurwitz systems, refer to  \cite{Auroux2003, GS, Matsu96, Moi81, ST03}, etc. 

Let $\iota$ be the mapping class of an involution of the fiber $f^{-1}(q_0)$ with 
$2g+2$ fixed points. 
The centralizer $HMG$ of $\iota$ in $MG$ 
is called the {\it hyperelliptic mapping class group} of $f^{-1}(q_0)$. 
A Lefschetz fibration is called {\it hyperelliptic} if the image of the monodromy representation $\rho$ is included in $HMG$.

\section{Main result} 

\begin{figure}[h]
\begin{center}
\mbox{\epsfxsize=8cm \epsfbox{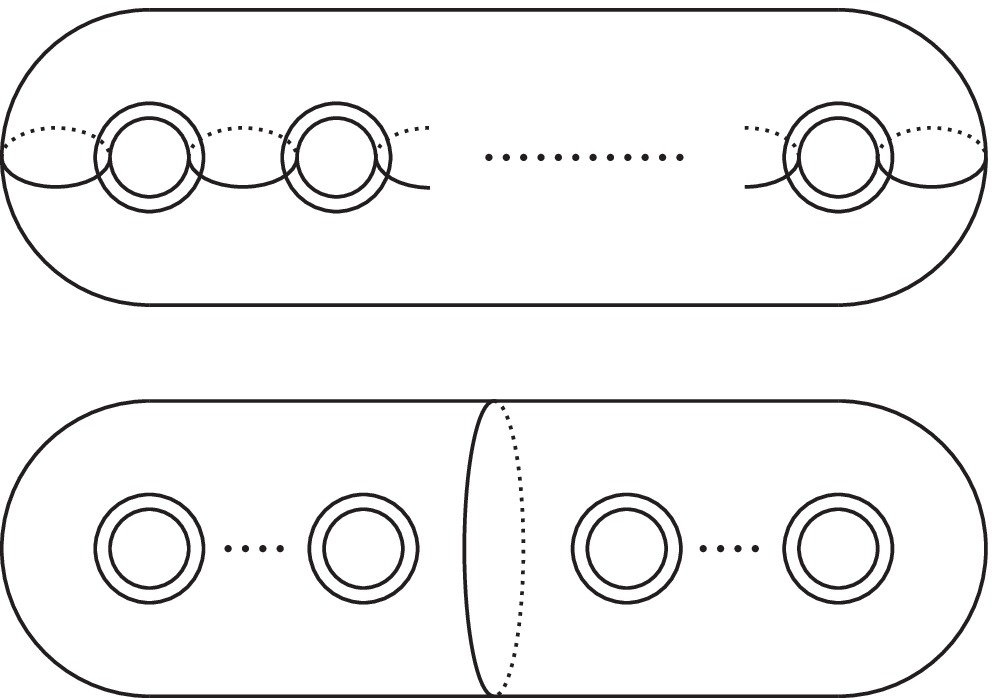}}
\put(-242,122){$C_1$}
\put(-200,140){$C_2$}
\put(-175,108){$C_3$}
\put(-155,140){$C_4$}
\put(-40,140){$C_{2g}$}
\put(2,122){$C_{2g+1}$}
\put(-200,50){$C_2$}
\put(-150,50){$C_{2h}$}
\put(-95,50){$C_{2h+2}$}
\put(-40,50){$C_{2g}$}
\put(-120,72){$S_h$}
\caption{\label{curves1} Curves on a general fiber}
\end{center}
\end{figure}

Let $\zeta_i$ $(i=1,\ldots, 2g+1)$ be positive Dehn twists along 
the loops $C_i$ $(i=1,\ldots, 2g+1)$ illustrated in Figure~\ref{curves1}.  
The hyperelliptic mapping class group $HMC$ of a genus-$g$ Riemann surface is generated by 
$\zeta_1, \ldots, \zeta_{2g+1}$, and 
the following relations are defining relations (cf. \cite{Birman}). 
{\allowdisplaybreaks %
\begin{eqnarray}
&& \zeta_i \zeta_j = \zeta_j \zeta_i     \quad \mbox{    if $|i-j|\geq 2$, } \label{eq:01} \\
&& \zeta_i \zeta_{i+1} \zeta_i = \zeta_{i+1} \zeta_i \zeta_{i+1}    \quad \mbox{    for $i=1, \dots, 2g$, } \label{eq:02} \\
&& \iota^2 =1    \quad \mbox{  where $\iota= \zeta_1 \cdots \zeta_{2g} \zeta_{2g+1}^2 \zeta_{2g} \cdots \zeta_1$,}\label{eq:03} \\
&&(\zeta_1 \cdots \zeta_{2g+1})^{2g+2}=1, \label{eq:04} \\
&&\iota  \,  \zeta_i =
 \zeta_i \,  \iota  
\quad \mbox{    for $i=1, \dots, 2g+1$.} \label{eq:05} 
\end{eqnarray}}
Let $\sigma_h$ be a positive Dehn twist along the loop $S_h$ illustrated in Figure~\ref{curves1}.  
Then $\sigma_h= (\zeta_1 \cdots \zeta_{2h})^{4h+2}$ for $h=1,\ldots ,[g/2]$. 

If $(g_1, \ldots, g_n)$ is a Hurwitz system of a genus-$g$ hyperelliptic Lefschetz fibration, then 
each $g_j$ is a conjugate of $\zeta_i$ or $\zeta_i^{-1}$, or a conjugate of $\sigma_h$ or $\sigma_h^{-1}$.  

Now we define basic Lefschetz fibrations.  

\begin{definition}[cf. \cite{Auroux2003, Auroux2005, Kam12, Matsu96, ST03}]{\rm 
{\it Basic Lefschetz fibrations},  $f_0$, $f_1$, $f_{2,h}$, $f'_1$ and $f'_{2,h}$,  are  
genus-$g$ hyperelliptic 
Lefschetz fibrations over $S^2$ whose Hurwitz systems are 
\begin{itemize}
\item[(1)] 
 $ W_0 = (T)^2 $ where $T = (\zeta_1, \zeta_2, \ldots, \zeta_{2g}, \zeta_{2g+1}, 
\zeta_{2g+1}, \zeta_{2g}, \ldots, \zeta_2, \zeta_1)$, 

\item[(2)]  $ W_1 = (\zeta_1, \zeta_2, \ldots, \zeta_{2g}, \zeta_{2g+1})^{2g+2}$, 

\item[(3)] $ W_{2,h} = (\zeta_{2g+1},\ldots ,\zeta_1, (\zeta_{2g-2h+1},\ldots ,\zeta_{2g+1}), \ldots , (\zeta_2, \ldots ,\zeta_{2h+2}), $ 

$(\zeta_1,\ldots ,\zeta_{2h+1}), \sigma_h, (\zeta_{2h+1},\ldots ,\zeta_1), $ 

$(\zeta_{2h+2},\ldots ,\zeta_2), \ldots , (\zeta_{2g+1},\ldots ,\zeta_{2g-2h+1}),  \zeta_1,\ldots ,\zeta_{2g+1}) $,  

\item[(4)]  $W'_1 = (\zeta_1, \zeta_1^{-1})$,  

\item[(5)]  $W'_{2,h} = (\sigma_h, \sigma_h^{-1})$, 
\end{itemize}
respectively.  
}\end{definition}

For example, $f_0$ has $4(2g+1)$ singular fibers, which are of type ${\rm I}^+$.  Thus $f_0$ is chiral and irreducible. 

\newcommand{\lw}[1]{\smash{\lower2.0ex\hbox{#1}}}
\begin{center}
\renewcommand{\arraystretch}{1.2}
 \begin{tabular}{|l@{\quad\vrule width0.8pt}c|c|c|c|c|c|}
 \hline 
 \lw{LF} & \multicolumn{4}{c|}{number of singular fibers} & \lw{chiral} & \lw{irreducible}   \\
 \cline{2-5}  & \, $n_0^+$ & $n_0^-$ & $n_k^+$ & $n_k^-$ & &\\
\noalign{\hrule height 0.8pt} 
 $f_0$ & $4(2g+1)$ & 0 & 0 & 0 & $\bigcirc$ & $\bigcirc$ \\
 $f_1$ & $\; 2(g+1)(2g+1)$ & 0 & 0 & 0 & $\bigcirc$ & $\bigcirc$ \\
 $f_{2,h}$ & $\; 8h(g-h)+4(2g+1)$ & 0 & $\delta_{hk}$ & 0 & $\bigcirc$ & $\times$ \\
 $f'_1$ & 1 & 1 & 0 &  0 & $\times$ & $\bigcirc$ \\
 $f'_{2,h}$ & 0 & 0 & $\delta_{hk}$ & $\delta_{hk}$ & $\times$ & $\times$ \\
\hline 
\end{tabular}
\end{center}

For two Lefschetz fibrations $f$ and $f'$ over $S^2$, we denote by $f \# f'$ a fiber-sum of $f$ and $f'$.  
By ${\#} m f$ for a positive integer $m$, we mean the fiber-sum of $m$ copies of $f$. 
If both $f$ and $f'$ are hyperelliptic, we assume that $f\# f'$ is also hyperelliptic. 

\begin{theorem}\label{thm:main}
Let $f$ be a genus-$g$ hyperelliptic Lefschetz fibration over $S^2$.  
Suppose that $n_h^+(f) \geq n_h^-(f)$ for $h=1,\ldots ,[g/2]$. Then 
\begin{itemize}
\item[{\rm (1)}] 
${\cal E}(f) := n_0^+(f) - n_0^-(f) -4\sum_{h=1}^{[g/2]}(n_h^+(f) - n_h^-(f))
(2h(g-h)+2g+1)$  
is a multiple of $2(2g+1)$ if $g$ is even, and that of $4(2g+1)$ if $g$ is odd.  
\end{itemize}  

\begin{itemize}
\item[{\rm (2)}] 
There exists a positive integer $m_0$ such that  
for any integer $m \geq m_0$, 
$$ 
f \,\# \, m   \,  f_0 
 \cong 
       {\#}  \, (a+m)  \,  f_0 
 \,  \# \,  b  \,  f_1 
 \,  \#  \, (\#_{h=1}^{[g/2]}c_h   \,  f_{2,h})
 \,  \#  \, d   \, f'_1
 \,  \#  \,  (\#_{h=1}^{[g/2]}e_h  \,  f'_{2,h}) $$
for some non-negative integers $a, b, c_1,\ldots ,c_{[g/2]}, d, e_1,\ldots ,e_{[g/2]}$. 
\item[{\rm (3)}] 
In $(2)$,  it holds that 
 $c_h = n_h^+(f) - n_h^-(f)$,  
$d = n_0^-(f)$ and $e_h = n_h^-(f)$.  
Although $a$ and $b$ are not determined uniquely,  we have 
$a = ({\cal E}(f) -2(g+1)(2g+1)b)/4(2g+1)$ and we can take 
$b \in \{0,1\}$. 
\end{itemize}
\end{theorem}

 \begin{remark}{\rm 
If $f$ is chiral, then 
$n_0^-(f) = n_1^-(f) \cdots =n_{[g/2]}^-(f)= 0$. 
By Theorem~\ref{thm:main}, we have \\ 
$$f    \,  \#   \,  m   \,  f_0  
\cong 
        \,  {\#}   \,  (a+m)   \,  f_0
  \,  {\#}   \,  b   \,  f_1 
  \,  {\#}   \,  c_1   \, f_{2,1} 
  \cdots 
  \,  {\#}  \,  c_{[g/2]}  \, f_{2,[g/2]} $$ 
for a sufficiently large integer $m$. 
Auroux and Smith \cite{AS} pointed out that a similar result follows from a work of 
Kharlamov and Kulikov \cite{KK}. 
}\end{remark}

 \begin{remark}{\rm 
If $g$ is even, then 
$b \equiv {\cal E}(f)/2(2g+1)\; ({\rm mod}\, 2)$. 
If $g$ is odd, the parity of $b$ is not determined by ${\cal E}(f)$. 
}\end{remark}


\section{Chart description} 

In this section we introduce a chart description for genus-$g$ hyperelliptic Lefschetz fibrations.  
We use the terminologies on chart description in \cite{Kam07}.  
For simplicity's sake, we only consider 
genus-$g$ hyperelliptic Lefschetz fibrations over $B$ such that $\partial B$ is empty or connected, and if  $\partial B$ is not empty, 
we assume that the monodromy along $\partial B$ is trivial.  Unless otherwise stated, genus-$g$ hyperelliptic Lefschetz fibrations over $B$ are assumed to be so.  

\begin{figure}[h]
\begin{center}
\mbox{\epsfxsize=10cm \epsfbox{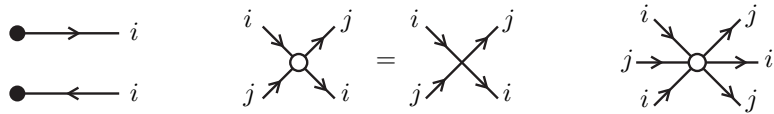}}
\put(-238,25){$i$}
\put(-238,2){$i$}
\put(-195,30){$i$}
\put(-158,30){$j$}
\put(-195,2){$j$}
\put(-158,2){$i$}
\put(-145,15){$=$}
\put(-132,30){$i$}
\put(-97,30){$j$}
\put(-132,2){$j$}
\put(-97,2){$i$}
\put(-45,32){$i$}
\put(-5,32){$j$}
\put(-52,15){$j$}
\put(2,15){$i$}
\put(-45,0){$i$}
\put(-5,0){$j$}
\end{center}
\vspace{-0.5cm}
\caption{Vertices of degree 1, 4, 6}
\label{vertices1}
\end{figure}

\begin{figure}[h]
\begin{center}
\mbox{\epsfxsize=10cm \epsfbox{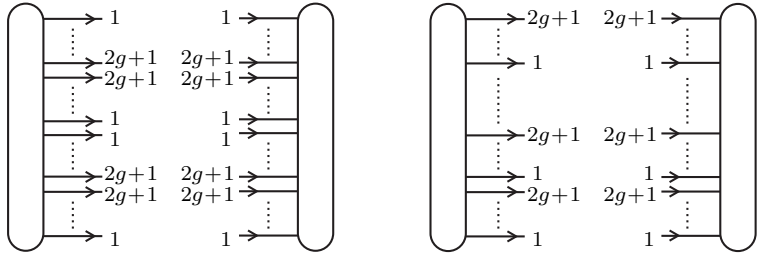}}
\put(-245,87){\footnotesize{$1$}}
\put(-247,72){\footnotesize{$2g\! +\! 1$}}
\put(-247,64){\footnotesize{$2g\! +\! 1$}}
\put(-245,49){\footnotesize{$1$}}
\put(-245,41){\footnotesize{$1$}}
\put(-247,28){\footnotesize{$2g\! +\! 1$}}
\put(-247,20){\footnotesize{$2g\! +\! 1$}}
\put(-245,3){\footnotesize{$1$}}
\put(-203,87){\footnotesize{$1$}}
\put(-218,72){\footnotesize{$2g\! +\! 1$}}
\put(-218,64){\footnotesize{$2g\! +\! 1$}}
\put(-203,49){\footnotesize{$1$}}
\put(-203,41){\footnotesize{$1$}}
\put(-218,28){\footnotesize{$2g\! +\! 1$}}
\put(-218,20){\footnotesize{$2g\! +\! 1$}}
\put(-203,3){\footnotesize{$1$}}
\put(-87,87){\footnotesize{$2g\! +\! 1$}}
\put(-85,70){\footnotesize{$1$}}
\put(-87,43){\footnotesize{$2g\! +\! 1$}}
\put(-85,28){\footnotesize{$1$}}
\put(-87,20){\footnotesize{$2g\! +\! 1$}}
\put(-85,3){\footnotesize{$1$}}
\put(-58,87){\footnotesize{$2g\! +\! 1$}}
\put(-43,70){\footnotesize{$1$}}
\put(-58,43){\footnotesize{$2g\! +\! 1$}}
\put(-43,28){\footnotesize{$1$}}
\put(-58,20){\footnotesize{$2g\! +\! 1$}}
\put(-43,3){\footnotesize{$1$}}
\end{center}
\vspace{-0.5cm}
\caption{Vertices of degree $4(2g+1)$, $2(g+1)(2g+1)$}
\label{vertices21}
\end{figure}

\begin{figure}[h]
\begin{center}
\mbox{\epsfxsize=10cm \epsfbox{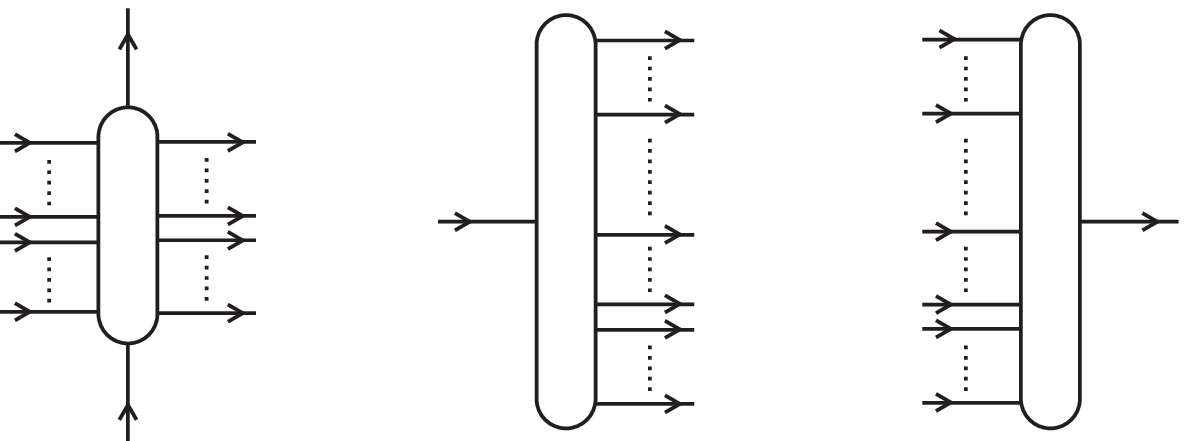}}
\put(-290,70){\footnotesize{$1$}}
\put(-305,53){\footnotesize{$2g\! +\! 1$}}
\put(-305,45){\footnotesize{$2g\! +\! 1$}}
\put(-290,28){\footnotesize{$1$}}
\put(-220,70){\footnotesize{$1$}}
\put(-220,53){\footnotesize{$2g\! +\! 1$}}
\put(-220,45){\footnotesize{$2g\! +\! 1$}}
\put(-220,28){\footnotesize{$1$}}
\put(-250,100){\footnotesize{$i$}}
\put(-250,0){\footnotesize{$i$}}
\put(-113,95){\footnotesize{$2h$}}
\put(-113,76){\footnotesize{$1$}}
\put(-113,48){\footnotesize{$2h$}}
\put(-113,31){\footnotesize{$1$}}
\put(-113,24){\footnotesize{$2h$}}
\put(-113,6){\footnotesize{$1$}}
\put(-75,95){\footnotesize{$2h$}}
\put(-70,76){\footnotesize{$1$}}
\put(-75,48){\footnotesize{$2h$}}
\put(-70,31){\footnotesize{$1$}}
\put(-75,24){\footnotesize{$2h$}}
\put(-70,6){\footnotesize{$1$}}
\put(-180,60){\footnotesize{$\sigma_h$}}
\put(-10,60){\footnotesize{$\sigma_h$}}
\end{center}
\vspace{-0.5cm}
\caption{Vertices of degree $2(4g+3)$, $4h(2h+1)+1$}
\label{vertices31}
\end{figure}

\begin{definition}[cf. \cite{Kam02, Kam07, Kam12, KMMW}]\label{def:chart} {\rm 
A {\it chart} in $B$ is a finite graph $\Gamma$ in $B$ (possibly being empty or having {\it hoops} that are closed edges without vertices) whose edges are labeled with an element of $\{1, \ldots , 2g+1, \sigma_1,\ldots ,\sigma_{[g/2]} \}$, and oriented so that the following conditions are satisfied (see Figure~\ref{vertices1}, \ref{vertices21}, \ref{vertices31}):   
\begin{itemize}

\item[(1)] The degree of each vertex is $1, 4, 6, 4(2g+1), 2(g+1)(2g+1), 2(4g+3)$ or $4h(2h+1)+1$.  
\item[(2)] For a degree-$1$ vertex, the adjacent edge is oriented outward or inward.  
\item[(3)] For a degree-$4$ vertex, two edges in each diagonal position have the same label and are oriented coherently; and the labels $i$ and $j$ of the diagonals are in $\{1, \dots, 2g+1\}$ with $| i - j | >1$.  
\item[(4)] For a degree-$6$ vertex, the six edges are alternately labeled $i$ and $j$  in $\{1, \dots, 2g+1\}$ with  $|i-j|=1$; and three consecutive edges are oriented outward while the other three are oriented inward.  
\item[(5)] For a vertex of degree $4(2g+1)$,  the edges are labeled with  $(1,\ldots ,2g+1,2g+1,\ldots ,1)^2$; and all edges are oriented outward or all edges are oriented inward. 
\item[(6)] For a vertex of degree $2(g+1)(2g+1)$, the edges are labeled with $(1,\ldots ,2g+1)^{2g+2}$ in a counterclockwise direction 
(or clockwise direction, resp.); and all edges are oriented outward (or inward, resp.). 
\item[(7)] For a vertex of degree $2(4g+3)$, the edges are labeled with $(1,\ldots ,2g+1,2g+1,\ldots ,1,i)^2$ in a counterclockwise direction where $i \in \{1, \dots, 2g+1\}$; and the first $4g+3$ edges are oriented outward and the latter ones are oriented inward.   
\item[(8)] For a vertex of degree $4h(2h+1)+1$, the edges are labeled with $((1,\ldots , 2h)^{4h+2}, \sigma_h)$  
in a counterclockwise direction 
(or clockwise direction, resp.); and the edges with labels $1,\ldots ,2h$ are oriented outward (or inward, resp.), and 
the edge with label $\sigma_h$ is oriented inward (or outward, resp.).  
\item[(9)] $\Gamma \cap \partial B = \emptyset$.   
\item[(10)] $\Gamma$ misses the base point $q_0 \in B$.  
\end{itemize}
}\end{definition}

\begin{remark}{\rm 
When we would treat genus-$g$ hyperelliptic Lefschetz fibrations over $B$ with $\partial B \neq \emptyset$ such that the monodromies along $\partial B$ are not  trivial, the condition (9) should be removed.  See \cite{Kam07}.  
}\end{remark}

We call a degree-$1$ vertex a {\it black vertex}.  
We say that a chart  is  {\it chiral} if every black vertex has an adjacent edge oriented outward.  
We say that a chart  is {\it irreducible} if there exist no edges with label $\sigma_h$.

For a chart $\Gamma$, let $\Delta_\Gamma$ be the set of black vertices.  
A chart $\Gamma$ determines a homomorphism $\pi_1(B \setminus \Delta_\Gamma, q_0) \to HMC$ as in \cite{Kam07}.  
By Theorem~5 of \cite{Kam07}, we have the following theorem. 

\begin{theorem}\label{thm:chartdescription}
Let $f $ be a genus-$g$ hyperelliptic Lefschetz fibration  over $B$, and let $\rho$ be the monodromy representation.    Then there is a chart $\Gamma$ in $B$ such that 
the monodromy representation $\rho$  equals 
the homomorphism  $\rho_\Gamma$ determined by $\Gamma$.    
\end{theorem}

A chart $\Gamma$ as in Theorem~\ref{thm:chartdescription} is called a {\it chart description} of $f$ or a {\it chart describing } $f$.  
Moreover, for a Hurwitz system $(g_1,\ldots ,g_n)$ of $f$ over $B=D^2$, we call $\Gamma$ a chart description of $(g_1,\ldots ,g_n)$. 
A chart $\Gamma$ in $D^2$ is also regarded as a chart in $S^2$ in the trivial way. 

\vspace{0.5cm}

We introduce some local moves on chart descriptions.  

(C1) For a chart $\Gamma$, suppose that there exists a chart $\Gamma'$ and an embedded 2-disk, say $E$,  in $B$ such that (i) $\partial E$ intersects with $\Gamma$ and $\Gamma'$ transversely  (or does not intersect with them) avoiding their vertices, (ii) $\Gamma$ and $\Gamma'$ have no black vertices in $E$, and (iii) $\Gamma$ and $\Gamma'$ are identical outside of $E$.  Then we say that $\Gamma'$ is obtained from $\Gamma$ by a C1-move.  

\begin{figure}[h]
\begin{center}
\mbox{\epsfxsize=9cm \epsfbox{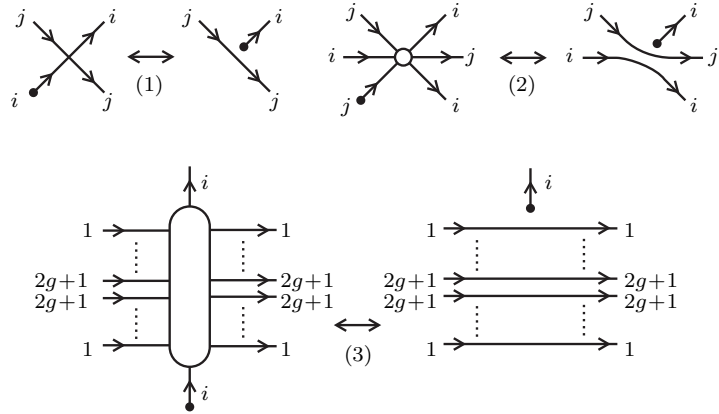}}
\put(-260,148){\footnotesize{$j$}}
\put(-262,116){\footnotesize{$i$}}
\put(-225,148){\footnotesize{$i$}}
\put(-227,116){\footnotesize{$j$}}
\put(-197,148){\footnotesize{$j$}}
\put(-162,148){\footnotesize{$i$}}
\put(-164,116){\footnotesize{$j$}}
\put(-215,122){\footnotesize{(1)}}
\put(-135,150){\footnotesize{$j$}}
\put(-95,150){\footnotesize{$i$}}
\put(-142,133){\footnotesize{$i$}}
\put(-90,133){\footnotesize{$j$}}
\put(-137,114){\footnotesize{$j$}}
\put(-95,114){\footnotesize{$i$}}
\put(-45,150){\footnotesize{$j$}}
\put(-52,133){\footnotesize{$i$}}
\put(-5,150){\footnotesize{$i$}}
\put(2,133){\footnotesize{$j$}}
\put(-5,114){\footnotesize{$i$}}
\put(-74,122){\footnotesize{(2)}}
\put(-235,67){\footnotesize{$1$}}
\put(-253,48){\footnotesize{$2g\! +\! 1$}}
\put(-253,40){\footnotesize{$2g\! +\! 1$}}
\put(-235,22){\footnotesize{$1$}}
\put(-160,67){\footnotesize{$1$}}
\put(-160,48){\footnotesize{$2g\! +\! 1$}}
\put(-160,40){\footnotesize{$2g\! +\! 1$}}
\put(-160,22){\footnotesize{$1$}}
\put(-190,85){\footnotesize{$i$}}
\put(-190,5){\footnotesize{$i$}}
\put(-105,67){\footnotesize{$1$}}
\put(-121,48){\footnotesize{$2g\! +\! 1$}}
\put(-121,40){\footnotesize{$2g\! +\! 1$}}
\put(-105,22){\footnotesize{$1$}}
\put(-30,67){\footnotesize{$1$}}
\put(-30,48){\footnotesize{$2g\! +\! 1$}}
\put(-30,40){\footnotesize{$2g\! +\! 1$}}
\put(-30,22){\footnotesize{$1$}}
\put(-60,85){\footnotesize{$i$}}
\put(-136,20){\footnotesize{(3)}}
\end{center}
\vspace{-0.5cm}
\caption{Some chart moves}
\label{moves1}
\end{figure}

(C2) For a chart, suppose that there is an edge $e$ joining a  degree-$4$ vertex and  a black vertex.  Remove the edge $e$ as in Figure~\ref{moves1}(1).  
We call this local move  a C2-move.  

(C3) For a chart, suppose that  there is an edge $e$ joining a degree-$6$ vertex and a black vertex.  Suppose that $e$ is neither the middle of three edges oriented outward nor the middle of the three edges oriented inward.  Then, remove the edge as in 
Figure~\ref{moves1}(2).  
We call this local move  a C3-move.  

(C4) In a chart, suppose that  there is an edge $e$ joining a degree-$2(4g+3)$ vertex and a black vertex.  Suppose that $e$ is one of the two edges labeled $i$ in Figure~\ref{vertices31}.  Then, remove the edge as in 
Figure~\ref{moves1}(3).  
We call this local move  a C4-move.  

\vspace{0.5cm}

When $\partial B \neq \emptyset$ and the base point $q_0$ is in $\partial B$, we introduce another move. 

(C5) Suppose that $\partial B \neq \emptyset$ and $q_0 \in \partial B$.  Let $\Gamma'$ be a chart that is the union of  a chart $\Gamma$ and some hoops which are parallel to and sufficiently near $\partial B$.  Then we say that $\Gamma'$ is obtained from $\Gamma$ by a C5-move.

\begin{definition}\label{def:chartmove}{\rm 
(1) {\it Chart moves} are C1-moves, C2-moves, C3-moves, C4-moves and their inverse moves.  

(2) Two charts in $B$ are said to be {\it chart move equivalent} (with respect to the base point $q_0$) if they are related by a finite sequence of chart moves and ambient isotopies of $B$ rel $q_0$, where we assume that chart moves are applied in embedded 2-disks in $B$ missing $q_0$. 

(3)   Two charts in $B$ are said to be {\it chart move equivalent up to conjugation} (with respect to the base point $q_0$) if they are related by a finite sequence of chart moves, C5-moves and ambient isotopies of $B$ rel $q_0$.  (It is not necessary to assume that chart moves are applied in embedded 2-disks in $B$ missing $q_0$.)
}\end{definition}

We say that two monodromy representations $\rho: \pi_1(B \setminus \Delta, q_0) \to HMC$ and $\rho': \pi_1(B \setminus \Delta', q_0) \to HMC$ are {\it equivalent} if there is a diffeomorphism $h : (B, q_0) \to (B, q_0)$ which is isotopic to the identity map rel $q_0$ such that $h(\Delta)= \Delta'$ and 
$\rho = \rho' \circ h_\ast$, where $h_\ast : \pi_1(B \setminus \Delta, q_0) \to  \pi_1(B \setminus \Delta', q_0)$ is the induced isomorphism.  

We say that two monodromy representations $\rho: \pi_1(B \setminus \Delta, q_0) \to HMC$ and $\rho': \pi_1(B \setminus \Delta', q_0) \to HMC$ are {\it equivalent up to conjugation} if there is an inner-automorphism of $HMC$, say $t$, and 
there is a diffeomorphism $h : (B, q_0) \to (B, q_0)$ which is isotopic to the identity map rel $q_0$ such that $h(\Delta)= \Delta'$ and 
$\rho = t \circ \rho' \circ h_\ast$.  

\vspace{0.5cm}

C1-moves in this paper are called  {\it chart moves of type $W$} in Definition~7 of \cite{Kam07}.  
C2-moves, C3-moves, C4-moves, C5-moves are not given explicitly in \cite{Kam07}.  
However, as shown in Fig. 22 and 23 of \cite{Kam07}, C2-moves and  C3-moves are 
equivalent to some 
local moves called {\it chart moves of transition} in Definition~14 of \cite{Kam07}.  
C4-moves are also equivalent to chart moves of transition in the sense of \cite{Kam07}.  
Thus, as stated in Section~8 of \cite{Kam07},  we see that if two charts are chart move equivalent in our sense (Definition~\ref{def:chartmove} (2)) then the monodromy representations determined by them are equivalent.  
C5-moves are equivalent to {\it chart moves of conjugacy} in (3) and (4) of Fig.~17 of \cite{Kam07}.  Again as in Section~7 of \cite{Kam07},  we see that if two charts are chart move equivalent up to conjugation (Definition~\ref{def:chartmove} (3)) then the monodromy representations determined by them are equivalent up to conjugation.  

Thus we have the following.  

\begin{theorem}
For two charts in $B$, if they are chart move equivalent (or chart move equivalent up to conjugation, resp.) then 
the monodromy representations determined by them are equivalent (or equivalent up to conjugation, resp.),  
and hence the hyperelliptic Lefschetz fibrations described by them are isomorphic.  
\end{theorem}

\begin{remark}{\rm 
By Theorem~16 of \cite{Kam07}, we see that two charts determine equivalent monodromy representations  if and only if they are related by C1-moves (chart move of type $W$), chart moves of transition, and ambient isotopies of $B$ rel $q_0$.  It is unknown to  the authors whether all chart moves of transition are consequence of our chart moves.  
}\end{remark}

We say that a black vertex of a chart $\Gamma$ is of {\it type ${\rm I}^+$}, 
{\it type ${\rm I}^-$},  {\it type ${\rm II}_h^+$} or  {\it type ${\rm II}_h^-$} if 
the adjacent edge is labeled in $\{1, \dots, 2g+1\}$ and oriented outward, if 
the adjacent edge is labeled in $\{1, \dots, 2g+1\}$ and oriented inward, if 
the adjacent edge is labeled $\sigma_h$ and oriented outward, or if 
the adjacent edge is labeled $\sigma_h$ and oriented inward, respectively.   

When $\Gamma$ is a chart description of a genus-$g$ hyperelliptic Lefschetz fibration $f : M \to B$,  
black vertices correspond to critical values of $f$, and 
the types of the vertices are the same with the types of the singular fibers over the corresponding critical values.   
For a chart $\Gamma$,  
we denote by 
$n_0^{+}(\Gamma)$, $n_0^{-}(\Gamma)$, $n_h^{+}(\Gamma)$, and $n_h^{-}(\Gamma)$, 
the numbers of black vertices of type ${\rm I}^+$, 
type ${\rm I}^-$,  type ${\rm II}_h^+$ and type ${\rm II}_h^-$, respectively.  They are equal to 
$n_0^{+}(f)$, $n_0^{-}(f)$, $n_h^{+}(f)$, and $n_h^{-}(f)$, respectively. 

\vspace{0.5cm}
If a chart $\Gamma$ is irreducible, then it is obvious that $n_h^{+}(\Gamma) = n_h^{-}(\Gamma) =0$ for $h=1,\ldots ,[g/2]$.  The converse is not true.  However we have the following.  

\begin{lemma}\label{lem:irreduciblechart}
Every chart $\Gamma$ with $n_h^{+}(\Gamma) = n_h^{-}(\Gamma) =0$ for $h=1,\ldots ,[g/2]$ is chart move equivalent to an irreducible chart.  
\end{lemma}

\begin{figure}[h]
\begin{center}
\mbox{\epsfxsize=12cm \epsfbox{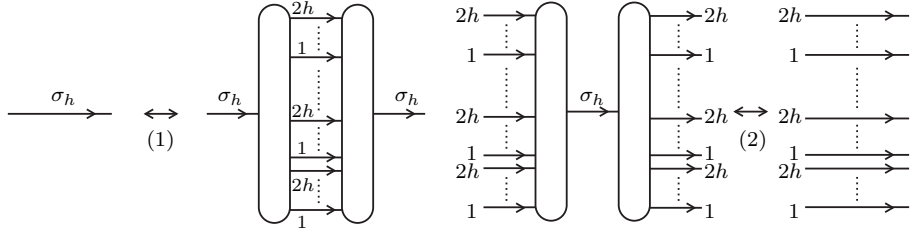}}
\put(-325,47){\footnotesize{$\sigma_h$}}
\put(-262,47){\footnotesize{$\sigma_h$}}
\put(-195,47){\footnotesize{$\sigma_h$}}
\put(-234,80){\scriptsize{$2h$}}
\put(-232,65){\scriptsize{$1$}}
\put(-234,41){\scriptsize{$2h$}}
\put(-232,27){\scriptsize{$1$}}
\put(-234,13){\scriptsize{$2h$}}
\put(-232,-1){\scriptsize{$1$}}
\put(-289,30){\footnotesize{(1)}}
\put(-125,47){\footnotesize{$\sigma_h$}}
\put(-172,77){\footnotesize{$2h$}}
\put(-168,62){\footnotesize{$1$}}
\put(-172,39){\footnotesize{$2h$}}
\put(-168,25){\footnotesize{$1$}}
\put(-172,18){\footnotesize{$2h$}}
\put(-168,3){\footnotesize{$1$}}
\put(-78,77){\footnotesize{$2h$}}
\put(-78,62){\footnotesize{$1$}}
\put(-78,39){\footnotesize{$2h$}}
\put(-78,25){\footnotesize{$1$}}
\put(-78,18){\footnotesize{$2h$}}
\put(-78,3){\footnotesize{$1$}}
\put(-50,77){\footnotesize{$2h$}}
\put(-46,62){\footnotesize{$1$}}
\put(-50,39){\footnotesize{$2h$}}
\put(-46,25){\footnotesize{$1$}}
\put(-50,18){\footnotesize{$2h$}}
\put(-46,3){\footnotesize{$1$}}
\put(-65,30){\footnotesize{(2)}}
\end{center}
\vspace{-0.5cm}
\caption{Chart moves}
\label{hoops1}
\end{figure}

{\it Proof.}  We can replace every hoop labeled $\sigma_h$ into $4h(2h+1)$ parallel hoops with labels $1,\ldots ,2h$ by 
a chart move depicted in Figure~\ref{hoops1} (1) followed by one  in Figure~\ref{hoops1} (2).  
Every edge labeled $\sigma_h$ whose endpoints are degree-$4h(2h+1)+1$ vertices is also removed by the latter move.  \qed  

\begin{proposition}\label{prop:chartdescriptionchiralirreducible}
A chiral (or irreducible, resp.) genus-$g$ hyperelliptic Lefschetz fibration has a chart description which is chiral 
(or irreducible, resp.).  
\end{proposition}

{\it Proof.} If $f$ is chiral, local monodromies around the critical values are all positive Dehn twists.  By the definition of a chart description, the adjacent edges of the black vertices are oriented outward.  Thus any chart description of $f$ is chiral.  If $f$ is 
irreducible, local monodromies around the critical values are Dehn twists along non-separating simple loops, which are conjugates of $\zeta_1, \dots, \zeta_{2g+1}$ and their inverses.  Thus any chart description $\Gamma$ of $f$ satisfies $n_h^{+}(\Gamma) = n_h^{-}(\Gamma) =0$ for $h=1,\ldots ,[g/2]$.  By Lemma~\ref{lem:irreduciblechart}, it changes to an irreducible one. \qed

\begin{figure}[h]
\begin{center}
\mbox{\epsfxsize=8cm \epsfbox{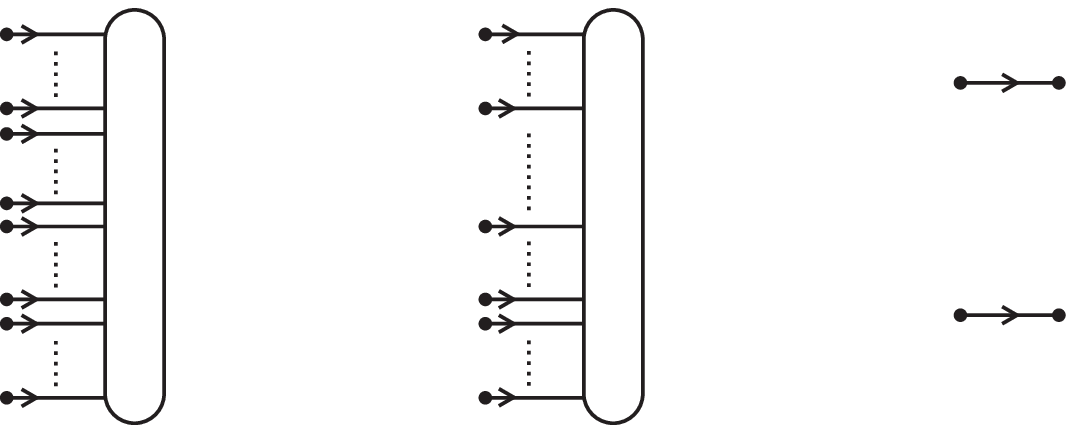}}
\put(-260,90){$N_0$}
\put(-233,81){\footnotesize{$1$}}
\put(-250,67){\footnotesize{$2g\! +\! 1$}}
\put(-250,59){\footnotesize{$2g\! +\! 1$}}
\put(-233,47){\footnotesize{$1$}}
\put(-233,39){\footnotesize{$1$}}
\put(-250,26){\footnotesize{$2g\! +\! 1$}}
\put(-250,18){\footnotesize{$2g\! +\! 1$}}
\put(-233,4){\footnotesize{$1$}}
\put(-210,-10){\footnotesize{deg$\, =4(2g\! +\! 1)$}}
\put(-170,90){$N_1$}
\put(-149,81){\footnotesize{$2g\! +\! 1$}}
\put(-133,65){\footnotesize{$1$}}
\put(-149,39){\footnotesize{$2g\! +\! 1$}}
\put(-133,26){\footnotesize{$1$}}
\put(-149,18){\footnotesize{$2g\! +\! 1$}}
\put(-133,4){\footnotesize{$1$}}
\put(-110,-10){\footnotesize{deg$\, =2(g\! +1)(2g\! +\! 1)$}}
\put(-50,85){$F_1$}
\put(-15,78){\footnotesize{$1$}}
\put(-50,35){$F_{2,h}$}
\put(-15,28){\footnotesize{$\sigma_h$}}
\end{center}
\vspace{-0.5cm}
\caption{Charts $N_0$, $N_1$, $F_1$ and $F_{2,h}$ describing 
$f_0$, $f_1$, $f'_1$ and $f'_{2,h}$}
\label{charts11}
\end{figure}

\begin{figure}[h]
\begin{center}
\mbox{\epsfxsize=12cm \epsfbox{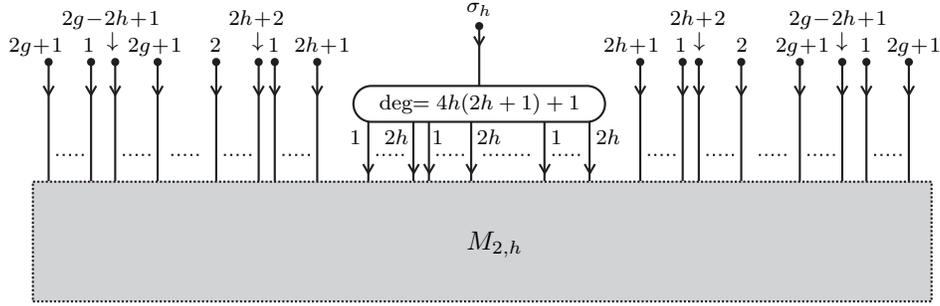}}
\put(-350,95){\footnotesize{$2g\! +\! 1$}}
\put(-322,95){\footnotesize{$1$}}
\put(-330,105){\footnotesize{$2g\! -\! 2h\! +\! 1$}}
\put(-313,97){\footnotesize{$\downarrow$}}
\put(-305,95){\footnotesize{$2g\! +\! 1$}}
\put(-274,95){\footnotesize{$2$}}
\put(-267,105){\footnotesize{$2h\! +\! 2$}}
\put(-258,97){\footnotesize{$\downarrow$}}
\put(-252,95){\footnotesize{$1$}}
\put(-242,95){\footnotesize{$2h\! +\! 1$}}
\put(-177,110){\footnotesize{$\sigma_h$}}
\put(-210,73){\footnotesize{deg$=4h(2h+1)+1$}}
\put(-177,20){$M_{2,h}$}
\put(-221,60){\footnotesize{$1$}}
\put(-208,60){\footnotesize{$2h$}}
\put(-190,60){\footnotesize{$1$}}
\put(-173,60){\footnotesize{$2h$}}
\put(-145,60){\footnotesize{$1$}}
\put(-128,60){\footnotesize{$2h$}}
\put(-125,95){\footnotesize{$2h\! +\! 1$}}
\put(-98,95){\footnotesize{$1$}}
\put(-100,105){\footnotesize{$2h\! +\! 2$}}
\put(-91,97){\footnotesize{$\downarrow$}}
\put(-75,95){\footnotesize{$2$}}
\put(-60,95){\footnotesize{$2g\! +\! 1$}}
\put(-55,105){\footnotesize{$2g\! -\! 2h\! +\! 1$}}
\put(-37,97){\footnotesize{$\downarrow$}}
\put(-28,95){\footnotesize{$1$}}
\put(-18,95){\footnotesize{$2g\! +\! 1$}}
\end{center}
\vspace{-0.5cm}
\caption{Chart $N_{2,h}$ describing $f_{2,h}$}
\label{charts21}
\end{figure}

In Figures~\ref{charts11} and \ref{charts21}, we show charts $N_0$, $N_1$, $N_{2,h}$, $F_1$ and $F_{2,h}$ describing 
$f_0$, $f_1$, $f_{2,h}$, $f'_1$ and $f'_{2,h}$.  
We call $N_0$ a (positive) {\it nucleon of degree-$4(2g+1)$} and $N_1$ a (positive) {\it nucleon of degree-$2(g+1)(2g+1)$}.  
The region named $M_{2,h}$ is a chart with some special properties (Lemma~\ref{lem:change}). 
A {\it free edge} means a chart consisting two black vertices and a single edge connecting them.  
$F_1$ and $F_{2,h}$ are free edges.  

\vspace{0.5cm}

Let $\Gamma$ and $\Gamma'$ be charts in $B=D^2$.  Divide $D^2$ into $2$-disks $D^2_1$ and $D^2_2$ by a properly embedded arc in $D^2$.  Put a small copy of $\Gamma$ in $D^2_1$ and a small copy of $\Gamma'$ in $D^2_2$.  We have a new chart in $D^2 = D^2_1 \cup D^2_2$.  We call it the {\it product} of $\Gamma$ and $\Gamma'$ and denote it by $\Gamma \oplus \Gamma'$.  
We say that $\Gamma$ is a {\it factor} of $\Gamma \oplus \Gamma'$.  
The chart $\Gamma \oplus \Gamma'$ is a chart description of the fiber sum $f \# f'$ of the Lefschetz fibrations $f$ and $f'$ described by $\Gamma$ and $\Gamma'$.   We denote by $n \Gamma$ the product $\Gamma \oplus \cdots \oplus \Gamma$ of $n$ copies of $\Gamma$.  (When $B=D^2$, the fiber sum $f \# f'$ of $f$ and $f'$ over $B$ is defined by using the boundary connected sum of the base spaces.)

\begin{theorem}\label{thm:chart1}
Let $\Gamma$ be a chart in $B=D^2$.  Suppose that $n_h^+(\Gamma) \geq n_h^-(\Gamma)$ for $h=1,\ldots ,[g/2]$. 
Then there exists a positive integer $m_0$ such that for any integer $m \geq m_0$, 
the chart $\Gamma \oplus m  \, N_0$ is chart move equivalent to 
$$\Gamma'  \oplus 
\left(\oplus_{h=1}^{[g/2]}(n_h^+(\Gamma)-  n_h^-(\Gamma)) N_{2,h}\right)  \oplus n_0^-(\Gamma) \, F_1  
\oplus  \left(\oplus_{h=1}^{[g/2]}n_h^-(\Gamma) F_{2,h}\right)  $$
for some chart $\Gamma'$ with $n_0^-(\Gamma')=n_h^+(\Gamma') = n_h^-(\Gamma')=0$ 
for $h=1,\ldots ,[g/2]$ such that $\Gamma'$ has $N_0$ as a factor.  
Moreover if $n_h^-(\Gamma)=0$ for $h=1,\ldots ,[g/2]$, we may take $m_0$ to be $n_0^-(\Gamma)  + \sum_{h=1}^{[g/2]} (h+1) n_h^+(\Gamma) +1 $.  
\end{theorem}

We prove Theorem~\ref{thm:chart1} in Section~\ref{sect:chart1proof}.  

\begin{corollary}\label{cor:chart1}
Let $f$ be a genus-$g$ hyperelliptic Lefschetz fibration over $B=D^2$ (or $S^2$) with $n_h^+(f) \geq n_h^-(f)$ for $h=1,\ldots ,[g/2]$. 
Then there exists a positive integer $m_0$ such that for any integer $m \geq m_0$, 
the fiber sum $f  \# m  f_0$ is equivalent to 
$$f'  \# 
\left(\#_{h=1}^{[g/2]}(n_h^+(f)-  n_h^-(f)) f_{2,h}\right)  \# n_0^-(f) \, f'_1  
\#  \left(\#_{h=1}^{[g/2]}n_h^-(f) f'_{2,h}\right)  $$ 
for some chiral and irreducible genus-$g$ hyperelliptic Lefschetz fibration $f'$ over $B=D^2$ (or $S^2$) 
such that the monodromy representation of $f'$ is transitive.  
Moreover if $n_h^-(f)=0$ for $h=1,\ldots ,[g/2]$, we may take $m_0$ to be $n_0^-(f)  + \sum_{h=1}^{[g/2]} (h+1) n_h^+(f) +1 $. 
\end{corollary}

\begin{figure}[h]
\begin{center}
\mbox{\epsfxsize=12.0cm \epsfbox{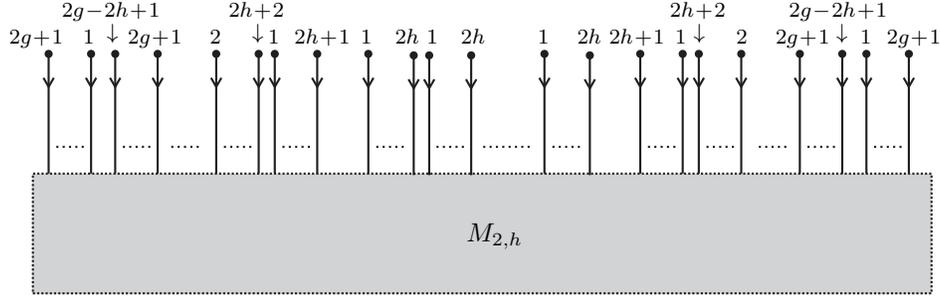}}
\put(-350,95){\footnotesize{$2g\! +\! 1$}}
\put(-322,95){\footnotesize{$1$}}
\put(-330,105){\footnotesize{$2g\! -\! 2h\! +\! 1$}}
\put(-313,97){\footnotesize{$\downarrow$}}
\put(-305,95){\footnotesize{$2g\! +\! 1$}}
\put(-274,95){\footnotesize{$2$}}
\put(-267,105){\footnotesize{$2h\! +\! 2$}}
\put(-258,97){\footnotesize{$\downarrow$}}
\put(-252,95){\footnotesize{$1$}}
\put(-242,95){\footnotesize{$2h\! +\! 1$}}
\put(-177,20){$M_{2,h}$}
\put(-217,95){\footnotesize{$1$}}
\put(-204,95){\footnotesize{$2h$}}
\put(-192,95){\footnotesize{$1$}}
\put(-179,95){\footnotesize{$2h$}}
\put(-150,95){\footnotesize{$1$}}
\put(-135,95){\footnotesize{$2h$}}
\put(-123,95){\footnotesize{$2h\! +\! 1$}}
\put(-98,95){\footnotesize{$1$}}
\put(-100,105){\footnotesize{$2h\! +\! 2$}}
\put(-91,97){\footnotesize{$\downarrow$}}
\put(-75,95){\footnotesize{$2$}}
\put(-60,95){\footnotesize{$2g\! +\! 1$}}
\put(-55,105){\footnotesize{$2g\! -\! 2h\! +\! 1$}}
\put(-37,97){\footnotesize{$\downarrow$}}
\put(-28,95){\footnotesize{$1$}}
\put(-18,95){\footnotesize{$2g\! +\! 1$}}
\end{center}
\vspace{-0.5cm}
\caption{Chart $P_{2,h}$, which is equivalent to $(h+1) N_0$}
\label{charts31}
\end{figure}

\begin{lemma}\label{lem:change}
There is a chart $M_{2,h}$ satisfying the following conditions. 
\begin{itemize}
\item[{\rm (1)}] 
$M_{2,h}$ consists of edges with labels in $\{1,\ldots ,2g+1\}$ and vertices whose degrees are in 
$\{4, 6, 4(2g+1), 2(4g+3)\}$. 
\item[{\rm (2)}] 
The chart $P_{2,h}$ depicted in Figure~$\ref{charts31}$ is chart move equivalent to 
$(h+1)N_0$. 
\end{itemize}
\end{lemma}

We prove Lemma~\ref{lem:change} in Section~\ref{sect:lemmaproof}.


\section{Proof of Theorem~\ref{thm:chart1}}\label{sect:chart1proof}

\begin{definition}{\rm 
A chart $\Gamma$ in a $2$-disk is {\it nomadic  with respect to} a chart $\Gamma_0$ in $B$ if 
for any two regions of the complement $B \setminus \Gamma_0$, say $R_1$ and $R_2$,  
the chart $\Gamma_0$ together with a small copy of $\Gamma$ in $R_1$ is chart   move equivalent to the chart $\Gamma_0$  together with a small copy of $\Gamma$ in $R_2$.  
A chart $\Gamma$ in a $2$-disk  is {\it nomadic} if it is nomadic with respect to every chart.  
}\end{definition}

\begin{figure}[h]
\begin{center}
\mbox{\epsfxsize=3cm \epsfbox{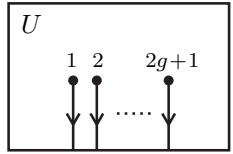}}
\put(-80,45){$U$}
\put(-63,32){\footnotesize{$1$}}
\put(-53,32){\footnotesize{$2$}}
\put(-33,32){\footnotesize{$2g\! +\! 1$}}
\end{center}
\vspace{-0.5cm}
\caption{Nomadic chart}
\label{nomadica}
\end{figure}

\begin{lemma}\label{lem:nomadic}
Let $D$ be a $2$-disk and $B$ a compact, connected and oriented surface.  
\begin{itemize}
\item[{\rm (1)}] 
Let $\Gamma$ be a chart in $D$.  
If there is a $2$-disk $U$ in $D$ such that $\Gamma \cap U$ is as in Figure~$\ref{nomadica}$, then
$\Gamma$ is nomadic.  
\item[{\rm (2)}] 
Let $\Gamma_0$ be a chart in $B$. 
If there is a $2$-disk $U$ in $B$ such that $\Gamma_0 \cap U$ is as in Figure~$\ref{nomadica}$, then 
any chart $\Gamma$ in a $2$-disk is nomadic with respect to $\Gamma_0$.  
\end{itemize}
\end{lemma}

\begin{figure}[h]
\begin{center}
\mbox{\epsfxsize=3cm \epsfbox{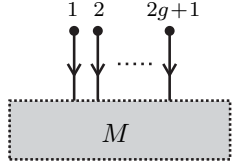}}
\put(-50,7){$M$}
\put(-63,55){\footnotesize{$1$}}
\put(-53,55){\footnotesize{$2$}}
\put(-33,55){\footnotesize{$2g\! +\! 1$}}
\end{center}
\vspace{-0.5cm}
\caption{Nomadic chart}
\label{nomadicb}
\end{figure}

\begin{figure}[h]
\begin{center}
\mbox{\epsfxsize=11cm \epsfbox{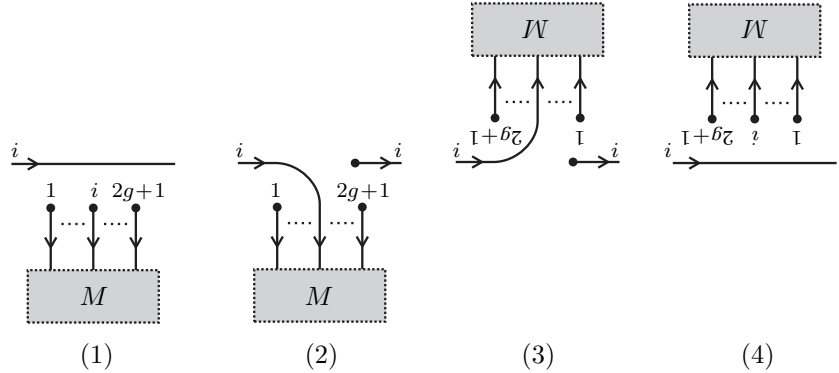}}
\put(-300,48){\footnotesize{$1$}}
\put(-283,48){\footnotesize{$i$}}
\put(-275,48){\footnotesize{$2g\! +\! 1$}}
\put(-313,63){\footnotesize{$i$}}
\put(-287,7){$M$}
\put(-287,-15){(1)}
\put(-215,48){\footnotesize{$1$}}
\put(-168,63){\footnotesize{$i$}}
\put(-190,48){\footnotesize{$2g\! +\! 1$}}
\put(-228,63){\footnotesize{$i$}}
\put(-202,7){$M$}
\put(-202,-15){(2)}
\put(-100,73){\rotatebox{180}{\footnotesize{$1$}}}
\put(-140,73){\rotatebox{180}{\footnotesize{$2g\! +\! 1$}}}
\put(-146,63){\footnotesize{$i$}}
\put(-86,63){\footnotesize{$i$}}
\put(-120,115){\rotatebox{180}{$M$}}
\put(-120,-15){(3)}
\put(-18,73){\rotatebox{180}{\footnotesize{$1$}}}
\put(-33,73){\rotatebox{180}{\footnotesize{$i$}}}
\put(-60,73){\rotatebox{180}{\footnotesize{$2g\! +\! 1$}}}
\put(-66,63){\footnotesize{$i$}}
\put(-37,115){\rotatebox{180}{$M$}}
\put(-37,-15){(4)}
\end{center}
\vspace{-0.5cm}
\caption{Chart moves}
\label{enter1}
\end{figure}

\begin{figure}[h]
\begin{center}
\mbox{\epsfxsize=9.0cm \epsfbox{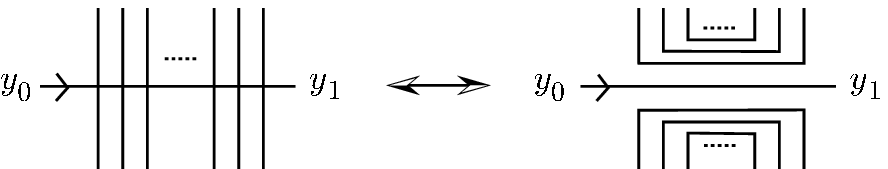}}
\end{center}
\vspace{-0.5cm}
\caption{Chart moves}
\label{sfg10}
\end{figure}

\begin{figure}[h]
\begin{center}
\mbox{\epsfxsize=9.0cm \epsfbox{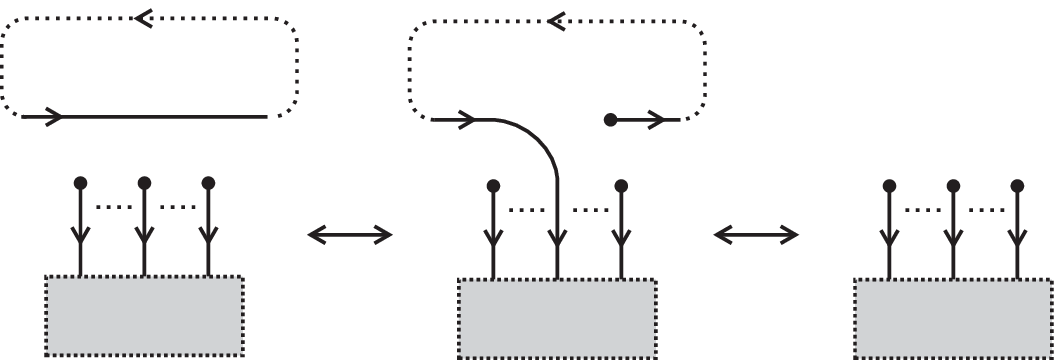}}
\put(-239,47){\footnotesize{$1$}}
\put(-223,47){\footnotesize{$i$}}
\put(-215,47){\footnotesize{$2g\! +\! 1$}}
\put(-233,62){\footnotesize{$i$}}
\put(-227,7){$M$}
\put(-139,47){\footnotesize{$1$}}
\put(-115,47){\footnotesize{$2g\! +\! 1$}}
\put(-133,62){\footnotesize{$i$}}
\put(-127,7){$M$}
\put(-42,47){\footnotesize{$1$}}
\put(-26,47){\footnotesize{$i$}}
\put(-18,47){\footnotesize{$2g\! +\! 1$}}
\put(-30,7){$M$}
\end{center}
\vspace{-0.5cm}
\caption{Chart moves}
\label{remove1}
\end{figure}

{\it Proof.} 
(1) First we consider a special case where $\Gamma$ is as in Figure~$\ref{nomadicb}$.  
Let $\Gamma_0$ be any chart in $B$, and put a small copy of $\Gamma$ in a region of $B \setminus \Gamma_0$.  
As shown in Figure~\ref{enter1}, it can pass through any edge of $\Gamma_0$ which is labeled in $\{1, \dots, 2g+1\}$.  
For an edge labeled $\sigma_h$, apply a chart move as in Figure~$\ref{hoops1} (1)$, let  $\Gamma$ pass through the $4(2h+1)$ edges with labels $1,\ldots ,2h$, and recover the edge labeled $\sigma_h$ by the move in Figure~\ref{hoops1}.   Thus we see that $\Gamma$ is nomadic.  
Now we consider a general case.  Take a point $y_0$ in the region $U$ and a point $y_1$ in the boundary $\partial D$.  Consider a simple path $\eta : [0,1] \to D$ connecting $y_0$ and $y_1$ such that $\eta$ intersects $\Gamma$ transversely.  Let $w$ be the intersection word of $\eta$ with respect to $\Gamma$ (see \cite{Kam02, Kam07}).  
Let $\Gamma'$ be a chart obtained from $\Gamma$ by adding some hoops surrounding $\Gamma$ such that 
the intersection word $w'$ of $\eta$ with respect to $\Gamma'$ is $w \cdot w^{-1}$.   Applying a chart move in a neighborhood of $\eta$ as in Figure~\ref{sfg10}, 
we have a chart $\Gamma''$ such that it coincides with $\Gamma'$ outside of the neighborhood of $\eta$ and the path $\eta$ misses $\Gamma''$. 
So $\Gamma''$ is as in  Figure~$\ref{nomadicb}$.  Note that $\Gamma''$ is chart move equivalent to $\Gamma$, since 
one can add or remove any hoop surrounding it by chart moves as in Figure~$\ref{remove1}$.  Since $\Gamma''$ is nomadic as shown in the previous case, we see that $\Gamma$ is nomadic.

\begin{figure}[h]
\begin{center}
\mbox{\epsfxsize=8.5cm \epsfbox{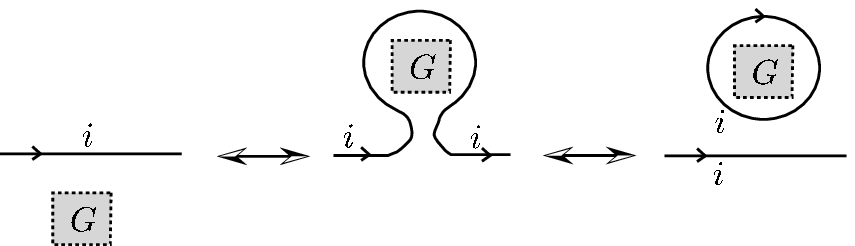}}
\end{center}
\vspace{-0.5cm}
\caption{Chart moves}
\label{sfg13}
\end{figure}

Now we prove (2). Let $U$ be a region such that $\Gamma_0 \cap U$ is as in Figure~$\ref{nomadica}$.  
It is sufficient to show that any chart $\Gamma$ put in a region of $B \setminus \Gamma_0$ can be moved into $U$.  
As shown in Figure~$\ref{sfg13}$,  $\Gamma$ can pass through any edge of $\Gamma_0$ by getting a surrounding hoop.   When $\Gamma$ arrives in $U$, it is surrounded some hoops, which can be removed by use of the edges of $\Gamma_0$ in $U$ as in Figure~$\ref{remove1}$.  \qed 

Now we prove Theorem~\ref{thm:chart1}.  

\begin{figure}[h]
\begin{center}
\mbox{\epsfxsize=7cm \epsfbox{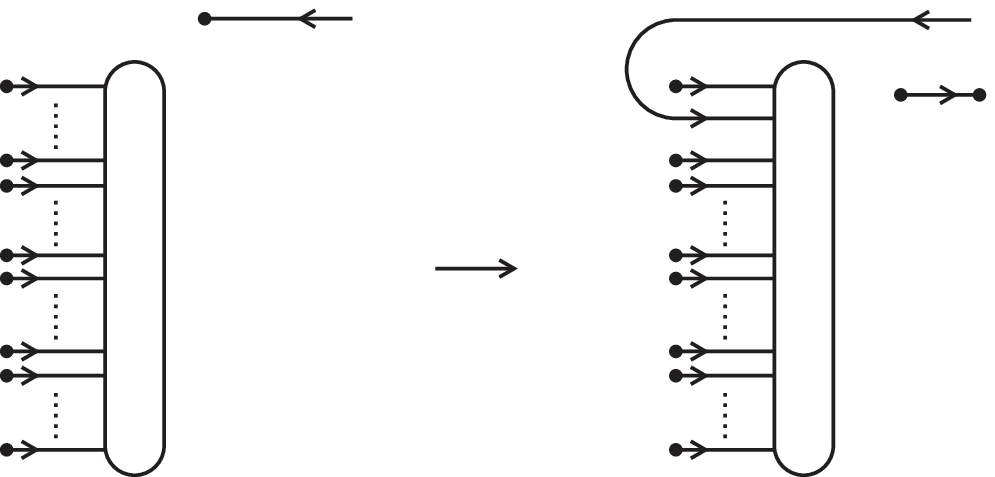}}
\put(-146,94){\footnotesize{$i$}}
\put(-206,77){\footnotesize{$1$}}
\put(-222,63){\footnotesize{$2g\! +\! 1$}}
\put(-222,55){\footnotesize{$2g\! +\! 1$}}
\put(-206,43){\footnotesize{$1$}}
\put(-206,36){\footnotesize{$1$}}
\put(-222,24){\footnotesize{$2g\! +\! 1$}}
\put(-222,16){\footnotesize{$2g\! +\! 1$}}
\put(-206,2){\footnotesize{$1$}}
\put(-162,13){\footnotesize{deg$=4(2g\! +\! 1)$}}
\put(-162,40){$N_0$}
\put(-11,94){\footnotesize{$i$}}
\put(-11,68){\footnotesize{$i$}}
\put(-70,79){\footnotesize{$1$}}
\put(-87,63){\footnotesize{$2g\! +\! 1$}}
\put(-87,55){\footnotesize{$2g\! +\! 1$}}
\put(-71,43){\footnotesize{$1$}}
\put(-71,36){\footnotesize{$1$}}
\put(-87,24){\footnotesize{$2g\! +\! 1$}}
\put(-87,16){\footnotesize{$2g\! +\! 1$}}
\put(-71,2){\footnotesize{$1$}}
\put(-27,13){\footnotesize{deg$=4(2g\! +\! 1)$}}
\end{center}
\vspace{-0.5cm}
\caption{Chart moves}
\label{moves31}
\end{figure}

\begin{figure}[h]
\begin{center}
\mbox{\epsfxsize=6.5cm \epsfbox{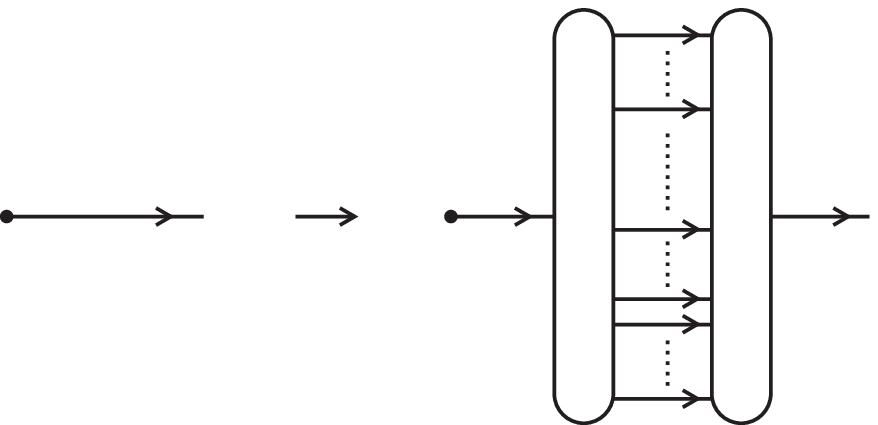}}
\put(-172,50){\footnotesize{$\sigma_h$}}
\put(-82,50){\footnotesize{$\sigma_h$}}
\put(-16,50){\footnotesize{$\sigma_h$}}
\put(-54,85){\scriptsize{$2h$}}
\put(-52,70){\scriptsize{$1$}}
\put(-54,44){\scriptsize{$2h$}}
\put(-52,29){\scriptsize{$1$}}
\put(-54,15){\scriptsize{$2h$}}
\put(-52,-1){\scriptsize{$1$}}
\put(-140,10){\footnotesize{deg$=4h(2h\! +\! 1)\! +\! 1$}}
\put(-18,10){\footnotesize{deg$=4h(2h\! +\! 1)\! +\! 1$}}
\end{center}
\vspace{-0.5cm}
\caption{Chart moves}
\label{moves41}
\end{figure}

{\it Proof of Theorem~\ref{thm:chart1}.}  
First we consider a case where $\Gamma$ is a chart with $n_h^-(\Gamma)=0$ for 
$h=1,\ldots ,[g/2]$.  
It suffices to show that 
$\Gamma \oplus (n_0^-(\Gamma)  + \sum_{h=1}^{[g/2]} (h+1) n_h^+(\Gamma) +1)   \, N_0$ is chart move equivalent to 
$$\Gamma'  \oplus 
\left(\oplus_{h=1}^{[g/2]}n_h^+(\Gamma) N_{2,h}\right)  \oplus n_0^-(\Gamma) \, F_1 $$
for some chart $\Gamma'$ with $n_0^-(\Gamma')=n_h^+(\Gamma') = n_h^-(\Gamma')=0$ 
such that $\Gamma'$ has $N_0$ as a factor.  
By Lemma~\ref{lem:nomadic}, $N_0$ is nomadic.  Thus we can move $N_0$ freely up to chart move equivalence.  
For each black vertex of type ${\rm I}^-$, move a chart $N_0$ near the vertex and apply a chart move as in Figure~\ref{moves31} to make a free edge.  Move the free edge toward the boundary of $B$ by the chart move as in Figure~\ref{sfg13}.  
Since there is at least one $N_0$ near $\partial B$, the hoops surrounding the free edge can be removed (Figure~\ref{remove1}), and we may also assume that the label of the free edge is $1$ (Lemma~18.24 of \cite{Kam02}).  
Thus we can change $\Gamma \oplus (n_0^-(\Gamma)  +  \sum_{h=1}^{[g/2]} (h+1) n_h^+(\Gamma) +1)   \, N_0$  so that 
all black vertices of type ${\rm I}^-$ are endpoints of $F_1$'s near $\partial B$.  
We still have $\sum_{h=1}^{[g/2]}(h+1)n_h^+(\Gamma) +1$ $N_0$'s near $\partial B$.  For each black vertex of type ${\rm II}_h^+$, move $h+1$ copies of $N_0$ near the vertex.   Change the copies of $N_0$ to a chart $P_{2,h}$ in Figure~$\ref{charts31}$ (Lemma~\ref{lem:change}).  
The edge adjacent to the vertex of type ${\rm II}_h^+$ is oriented outward and is labeled $\sigma_h$.  Apply a chart move as in Figure~\ref{moves41}, and then apply a chart move between the $4(2h+1)$ edges there and the $4(2h+1)$ edges of $P_{2,h}$ to get one $N_{2,h}$.  Move the chart $N_{2,h}$ toward $\partial B$.  (Note that $N_{2,h}$ is nomadic by Lemma~\ref{lem:nomadic}.)  Now all black vertices of type ${\rm II}_h^+$ belong to $N_{2,h}$'s near $\partial B$.  
We still have one $N_0$ near $\partial B$.  Thus the chart is 
$\Gamma'  \oplus 
(\oplus_{h=1}^{[g/2]}n_h^+(\Gamma) N_{2,h})  \oplus n_0^-(\Gamma) \, F_1 $ for a chart $\Gamma'$ with $n_0^-(\Gamma')=n_h^+(\Gamma') = n_h^-(\Gamma')=0$ 
such that $\Gamma'$ has $N_0$ as a factor.  

We consider a case where $\Gamma$ is a chart with $n_h^+(\Gamma) \geq n_h^-(\Gamma)>0$ for $h=1,\ldots ,[g/2]$.  
Let $v$ be a black vertex of type ${\rm II}_h^-$.  Choose a black vertex $v'$ of type ${\rm II}_h^+$ and consider a simple path $\eta$ from $v$ to $v'$.  If $\eta$ intersects an edge labeled $\sigma_h$, then apply a chart move depicted in Figure~$\ref{hoops1} (1)$ and we assume that $\eta$ intersects only edges with labels in $\{1, \dots, 2g+1\}$.  For each intersection of $\eta$ and the chart, we assert one $N_0$ and apply a chart move as in Figure~$\ref{enter1} (2)$ so that $\eta$ does not intersect the chart.  Now move $v$ along $\eta$ toward $v'$ and by a chart move we can make a free edge with label $\sigma_h$, that is $F_{2,h}$.  Move this $F_{2,h}$ toward $\partial B$ by 
moves as in Figure~\ref{sfg13}.   The hoops surrounding the free edge can be removed by adding one $N_0$ near $\partial B$ as before.  By this procedure, we can move all black vertices of type ${\rm II}_h^-$ near $\partial B$ as endpoints of $F_{2,h}$'s.  
The number of $F_{2,h}$'s is $n_h^-(\Gamma)$.  
There are $n_h^+(\Gamma) - n_h^-(\Gamma)$ 
black vertices of type ${\rm II}_h^+$ in the chart, besides the endpoints of $F_{2,h}$'s.  For each black vertex of type ${\rm II}_h^+$, that is not an endpoint of $F_{2,h}$, add $h+1$ copies of $N_0$ to make $P_{2,h}$.  As in the previous case, we can move the black vertex of type ${\rm II}_h^+$ as an endpoint of $N_{2,h}$ near $\partial B$.  The number of $N_{2,h}$'s is $n_h^+(\Gamma) - n_h^-(\Gamma)$.  As in the previous case, we move black vertices of type ${\rm I}^-$ as endpoints of $F_1$'as near $\partial B$.  The number of $F_1$'s is $n_0^-(\Gamma)$.   Thus we have a chart written as 
$$\Gamma'  \oplus 
\left(\oplus_{h=1}^{[g/2]}(n_h^+(\Gamma)-  n_h^-(\Gamma)) N_{2,h}\right)  \oplus n_0^-(\Gamma) \, F_1  
\oplus  \left(\oplus_{h=1}^{[g/2]}n_h^-(\Gamma) F_{2,h}\right)  $$
for some chart $\Gamma'$ with $n_0^-(\Gamma')=n_h^+(\Gamma') = n_h^-(\Gamma')=0$ 
such that $\Gamma'$ has $N_0$ as a factor.  \qed


\section{Proof of Theorem~\ref{thm:main}}\label{sect:mainproof}

\begin{proposition}\label{prop:chart2}
Let $\Gamma$ be a chart description of a chiral and irreducible genus-$g$ hyperelliptic  Lefschetz fibration over $B=D^2$ (or $S^2$).  
There exists a positive integer $m$ such that 
$\Gamma \oplus m N_0$ is chart move equivalent to $ (a + m) N_0 \oplus b N_1$ for some 
integers $a$ and $b$.  
\end{proposition}

\begin{figure}[h]
\begin{center}
\mbox{\epsfxsize=8cm \epsfbox{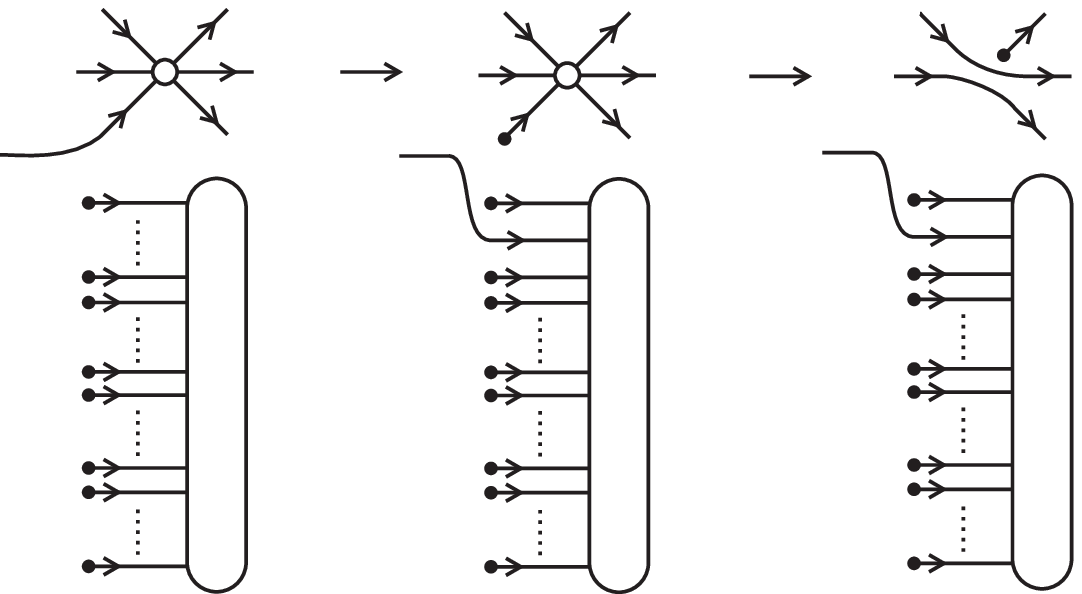}}
\put(-245,120){\footnotesize{deg$=6$}}
\put(-211,123){\footnotesize{$j$}}
\put(-216,109){\footnotesize{$i$}}
\put(-226,85){\footnotesize{$j$}}
\put(-178,123){\footnotesize{$i$}}
\put(-172,109){\footnotesize{$j$}}
\put(-178,92){\footnotesize{$i$}}
\put(-216,80){\footnotesize{$1$}}
\put(-232,66){\footnotesize{$2g\! +\! 1$}}
\put(-232,58){\footnotesize{$2g\! +\! 1$}}
\put(-216,46){\footnotesize{$1$}}
\put(-216,39){\footnotesize{$1$}}
\put(-232,26){\footnotesize{$2g\! +\! 1$}}
\put(-232,18){\footnotesize{$2g\! +\! 1$}}
\put(-216,3){\footnotesize{$1$}}
\put(-172,40){$N_0$}
\put(-126,123){\footnotesize{$j$}}
\put(-131,109){\footnotesize{$i$}}
\put(-128,97){\footnotesize{$j$}}
\put(-141,85){\footnotesize{$j$}}
\put(-93,123){\footnotesize{$i$}}
\put(-87,109){\footnotesize{$j$}}
\put(-93,92){\footnotesize{$i$}}
\put(-111,85){\footnotesize{$1$}}
\put(-147,66){\footnotesize{$2g\! +\! 1$}}
\put(-147,58){\footnotesize{$2g\! +\! 1$}}
\put(-131,46){\footnotesize{$1$}}
\put(-131,39){\footnotesize{$1$}}
\put(-147,26){\footnotesize{$2g\! +\! 1$}}
\put(-147,18){\footnotesize{$2g\! +\! 1$}}
\put(-131,3){\footnotesize{$1$}}
\put(-38,123){\footnotesize{$j$}}
\put(-43,109){\footnotesize{$i$}}
\put(-53,86){\footnotesize{$j$}}
\put(-5,123){\footnotesize{$i$}}
\put(1,109){\footnotesize{$j$}}
\put(-5,92){\footnotesize{$i$}}
\put(-23,86){\footnotesize{$1$}}
\put(-59,67){\footnotesize{$2g\! +\! 1$}}
\put(-59,59){\footnotesize{$2g\! +\! 1$}}
\put(-43,47){\footnotesize{$1$}}
\put(-43,40){\footnotesize{$1$}}
\put(-59,27){\footnotesize{$2g\! +\! 1$}}
\put(-59,19){\footnotesize{$2g\! +\! 1$}}
\put(-43,4){\footnotesize{$1$}}
\end{center}
\vspace{-0.5cm}
\caption{Chart moves}
\label{moves51}
\end{figure}

\begin{figure}[h]
\begin{center}
\mbox{\epsfxsize=9.5cm \epsfbox{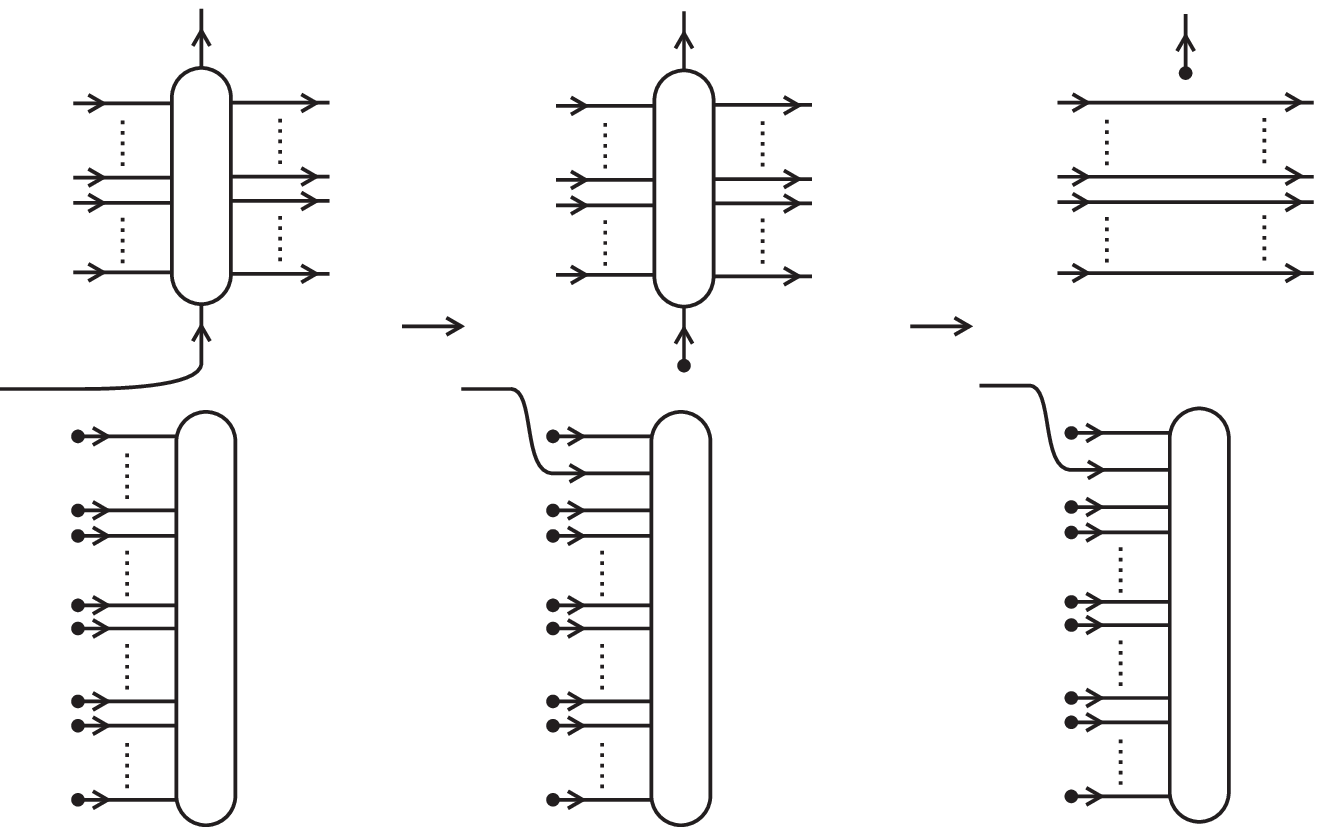}}
\put(-295,170){\footnotesize{deg$=2(4g\! +\! 3)$}}
\put(-261,147){\footnotesize{$1$}}
\put(-279,133){\footnotesize{$2g\! +\! 1$}}
\put(-279,125){\footnotesize{$2g\! +\! 1$}}
\put(-261,112){\footnotesize{$1$}}
\put(-200,147){\footnotesize{$1$}}
\put(-200,133){\footnotesize{$2g\! +\! 1$}}
\put(-200,125){\footnotesize{$2g\! +\! 1$}}
\put(-200,112){\footnotesize{$1$}}
\put(-237,160){\footnotesize{$j$}}
\put(-271,83){\footnotesize{$j$}}
\put(-262,78){\footnotesize{$1$}}
\put(-278,65){\footnotesize{$2g\! +\! 1$}}
\put(-278,57){\footnotesize{$2g\! +\! 1$}}
\put(-262,45){\footnotesize{$1$}}
\put(-262,38){\footnotesize{$1$}}
\put(-278,25){\footnotesize{$2g\! +\! 1$}}
\put(-278,17){\footnotesize{$2g\! +\! 1$}}
\put(-262,3){\footnotesize{$1$}}
\put(-217,40){$N_0$}
\put(-161,147){\footnotesize{$1$}}
\put(-178,133){\footnotesize{$2g\! +\! 1$}}
\put(-178,125){\footnotesize{$2g\! +\! 1$}}
\put(-160,112){\footnotesize{$1$}}
\put(-100,147){\footnotesize{$1$}}
\put(-100,133){\footnotesize{$2g\! +\! 1$}}
\put(-100,125){\footnotesize{$2g\! +\! 1$}}
\put(-100,112){\footnotesize{$1$}}
\put(-137,160){\footnotesize{$j$}}
\put(-137,97){\footnotesize{$j$}}
\put(-176,83){\footnotesize{$j$}}
\put(-145,83){\footnotesize{$1$}}
\put(-181,65){\footnotesize{$2g\! +\! 1$}}
\put(-181,57){\footnotesize{$2g\! +\! 1$}}
\put(-165,45){\footnotesize{$1$}}
\put(-165,38){\footnotesize{$1$}}
\put(-181,25){\footnotesize{$2g\! +\! 1$}}
\put(-181,17){\footnotesize{$2g\! +\! 1$}}
\put(-165,3){\footnotesize{$1$}}
\put(-60,147){\footnotesize{$1$}}
\put(-76,133){\footnotesize{$2g\! +\! 1$}}
\put(-76,125){\footnotesize{$2g\! +\! 1$}}
\put(-60,112){\footnotesize{$1$}}
\put(2,147){\footnotesize{$1$}}
\put(2,133){\footnotesize{$2g\! +\! 1$}}
\put(2,125){\footnotesize{$2g\! +\! 1$}}
\put(2,112){\footnotesize{$1$}}
\put(-35,160){\footnotesize{$j$}}
\put(-70,83){\footnotesize{$j$}}
\put(-38,84){\footnotesize{$1$}}
\put(-74,65){\footnotesize{$2g\! +\! 1$}}
\put(-74,57){\footnotesize{$2g\! +\! 1$}}
\put(-58,45){\footnotesize{$1$}}
\put(-58,38){\footnotesize{$1$}}
\put(-74,25){\footnotesize{$2g\! +\! 1$}}
\put(-74,17){\footnotesize{$2g\! +\! 1$}}
\put(-58,4){\footnotesize{$1$}}
\end{center}
\vspace{-0.5cm}
\caption{Chart moves}
\label{moves61}
\end{figure}

\begin{figure}[h]
\begin{center}
\mbox{\epsfxsize=11cm \epsfbox{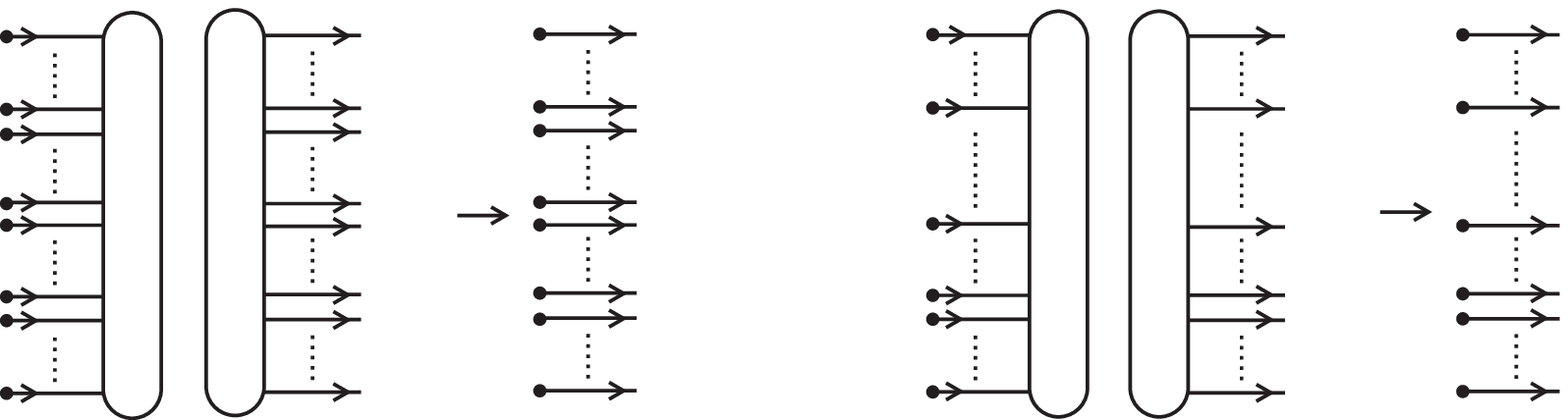}}
\put(-318,75){\footnotesize{$1$}}
\put(-334,62){\footnotesize{$2g\! +\! 1$}}
\put(-334,55){\footnotesize{$2g\! +\! 1$}}
\put(-318,43){\footnotesize{$1$}}
\put(-318,36){\footnotesize{$1$}}
\put(-334,24){\footnotesize{$2g\! +\! 1$}}
\put(-334,17){\footnotesize{$2g\! +\! 1$}}
\put(-318,3){\footnotesize{$1$}}
\put(-240,75){\footnotesize{$1$}}
\put(-240,62){\footnotesize{$2g\! +\! 1$}}
\put(-240,55){\footnotesize{$2g\! +\! 1$}}
\put(-240,43){\footnotesize{$1$}}
\put(-240,36){\footnotesize{$1$}}
\put(-240,24){\footnotesize{$2g\! +\! 1$}}
\put(-240,17){\footnotesize{$2g\! +\! 1$}}
\put(-240,3){\footnotesize{$1$}}
\put(-185,75){\footnotesize{$1$}}
\put(-185,62){\footnotesize{$2g\! +\! 1$}}
\put(-185,55){\footnotesize{$2g\! +\! 1$}}
\put(-185,43){\footnotesize{$1$}}
\put(-185,36){\footnotesize{$1$}}
\put(-185,24){\footnotesize{$2g\! +\! 1$}}
\put(-185,17){\footnotesize{$2g\! +\! 1$}}
\put(-185,3){\footnotesize{$1$}}
\put(-300,-10){$N_0$}
\put(-270,-10){\footnotesize{deg$=4(2g\! +\! 1)$}}
\put(-148,75){\footnotesize{$2g\! +\! 1$}}
\put(-132,60){\footnotesize{$1$}}
\put(-148,37){\footnotesize{$2g\! +\! 1$}}
\put(-132,24){\footnotesize{$1$}}
\put(-148,17){\footnotesize{$2g\! +\! 1$}}
\put(-132,3){\footnotesize{$1$}}
\put(-55,75){\footnotesize{$2g\! +\! 1$}}
\put(-55,60){\footnotesize{$1$}}
\put(-55,35){\footnotesize{$2g\! +\! 1$}}
\put(-55,24){\footnotesize{$1$}}
\put(-55,17){\footnotesize{$2g\! +\! 1$}}
\put(-55,3){\footnotesize{$1$}}
\put(2,75){\footnotesize{$2g\! +\! 1$}}
\put(2,60){\footnotesize{$1$}}
\put(2,37){\footnotesize{$2g\! +\! 1$}}
\put(2,24){\footnotesize{$1$}}
\put(2,17){\footnotesize{$2g\! +\! 1$}}
\put(2,3){\footnotesize{$1$}}
\put(-115,-10){$N_1$}
\put(-90,-10){\footnotesize{deg$=4(g\! +\! 1)(2g\! +\! 1)$}}
\end{center}
\vspace{-0.5cm}
\caption{Chart moves}
\label{moves71}
\end{figure}

{\it Proof.}  Since $f$ is chiral and irreducible,  we may assume that $\Gamma$ is chiral and irreducible by 
Proposition~\ref{prop:chartdescriptionchiralirreducible}.  
Adding some $N_0$'s to the chart and applying chart moves shown 
in Figures~\ref{moves51}--\ref{moves71},  
we can remove all degree-$6$ vertices, degree-$2(4g+3)$ vertices,  degree-$4(2g+1)$ vertices whose adjacent edges are oriented outward, and 
degree-$2(g+1)(2g+1)$ vertices whose adjacent edges are oriented outward.    
(Since it is easily seen that $(g+1) N_0$ is chart move equivalent to  $2 N_1$ (cf. \cite{Endo}), we may add $N_1$'s too.)  
Remove all hoops using an $N_0$ (Figure~$\ref{remove1}$).  
Now every edge is adjacent to a black vertex, a degree-$4$ vertex, a degree-$4(2g+1)$ vertex whose adjacent edges are oriented inward or a degree-$2(g+1)(2g+1)$ vertex whose adjacent edges are oriented inward.   Note that for a degree-$4$ vertex, the two incoming adjacent edges have black vertices at the other end.  Thus by a C2-move (Figure~\ref{moves1}(1)), we can remove the degree-$4$ vertex.  Remove all degree-$4$ vertices this way.  Now the chart is a union of some $N_0$'s and $N_1$'s.  \qed

Now we have a corollary to Proposition~\ref{prop:chart2}.  

\begin{corollary}\label{cor:chart2}
Let $f$ be a chiral and irreducible genus-$g$ hyperelliptic Lefschetz fibration over $S^2$.  
There exists a positive number $m$ such that 
$f \#  m f_0 \cong (a + m) f_0 \# b f_1$ for some 
integers $a$ and $b$.  
\end{corollary} 

Using this corollary, we have a proof of Theorem~\ref{thm:main}.

{\it Proof of Theorem~\ref{thm:main}.}  
Corollary~\ref{cor:chart1} and Corollary~\ref{cor:chart2} imply the assertions (2) of Theorem~\ref{thm:main}.    
We shall compare the number of singular fibers of each type of $f \,\# \, m   \,  f_0 $ with that of 
$ 
{\#}  \, (a+m)  \,  f_0 
 \,  \# \,  b  \,  f_1 
 \,  \#  \, (\#_{h=1}^{[g/2]}c_h   \,  f_{2,h})
 \,  \#  \, d   \, f'_1
 \,  \#  \,  (\#_{h=1}^{[g/2]}e_h  \,  f'_{2,h}) $. 
We have already used the information on the numbers of singular fibers of type ${\rm I}^-$, ${\rm II}_h^+$ and type ${\rm II}_h^-$ to determine $c_h$, $d$ and $e_h$; $c = n_h^+(f) - n_h^-(f)$,  
$d = n_0^-(f)$ and $e = n_h^-(f)$.   
The number of singular fibers of type ${\rm I}^+$ of 
$f \,\# \, m   \,  f_0 $ is $n_0^+(f) + 4(2g+1) m$, and that of 
$
{\#}  \, (a+m)  \,  f_0 
 \,  \# \,  b  \,  f_1 
 \,  \#  \, (\#_{h=1}^{[g/2]}c_h   \,  f_{2,h})
 \,  \#  \, d   \, f'_1
 \,  \#  \,  (\#_{h=1}^{[g/2]}e_h  \,  f'_{2,h}) $
is $4(2g+1) (a+m) + 2(g+1)(2g+1) b + \sum_{h=1}^{[g/2]}(8h(g-h)+4(2g+1))c_h + d$.  
From this equality, we have 
\begin{align*}
& \; 4(2g+1) a + 2(g+1)(2g+1) b \\
= & \;  n_0^+(f) - n_0^-(f) 
-4\sum_{h=1}^{[g/2]}(n_h^+(f)-n_h^-(f))(2h(g-h)+2g+1).
\end{align*}
Thus the right hand side, which is ${\cal E}(f)$, is a multiple of $2(2g+1)$ if $g$ is even, 
and that of $4(2g+1)$ if $g$ is odd.  
It is well-known that $(g+1) f_0 \cong 2 f_1$ (cf. \cite{Endo}).   
Therefore we have $a = ({\cal E}(f) -2(g+1)(2g+1) b)/4(2g+1)$ and we can take 
$b$ to be $0$ or $1$.   \qed


\section{Proof of Lemma~\ref{lem:change}}\label{sect:lemmaproof}

We introduce some local moves on Hurwitz systems.

(H1) For a Hurwitz system, suppose that there are two consecutive components 
$(\zeta_i,\zeta_j)$ with $|i-j|>1$. Replace them by $(\zeta_j,\zeta_i)$. 
We call this local move  an H1-move.  

(H2) For a Hurwitz system, suppose that there are three consecutive components 
$(\zeta_i,\zeta_j,\zeta_i)$ with $|i-j|=1$. Replace them by $(\zeta_j,\zeta_i,\zeta_j)$. 
We call this local move  an H2-move.  

(H3) For a Hurwitz system, suppose that there are $4g+3$ consecutive components 
$(\zeta_1,\ldots ,\zeta_{2g+1},\zeta_{2g+1},\ldots ,\zeta_1,\zeta_i)$ where $i \in \{1, \dots, 2g+1\}$. 
Replace them by $(\zeta_i,\zeta_1,\ldots ,\zeta_{2g+1},\zeta_{2g+1},\ldots ,\zeta_1)$. 
We call this local move  an H3-move.

\begin{lemma}\label{lem:hmoves}
Let $\Gamma$ be a chart description of a Hurwitz system $(g_1, \ldots, g_n)$. 
If a Hurwitz system $(g'_1, \ldots, g'_n)$ is obtained from $(g_1, \ldots, g_n)$ by a H1-move (or H2-move, H3-move, resp.) or its inverse, a chart description $\Gamma'$ is obtained from 
$\Gamma$ by a C2-move (or C3-move, C4-move, resp.) or its inverse. 
\end{lemma}

\begin{figure}[h]
\begin{center}
\mbox{\epsfxsize=9cm \epsfbox{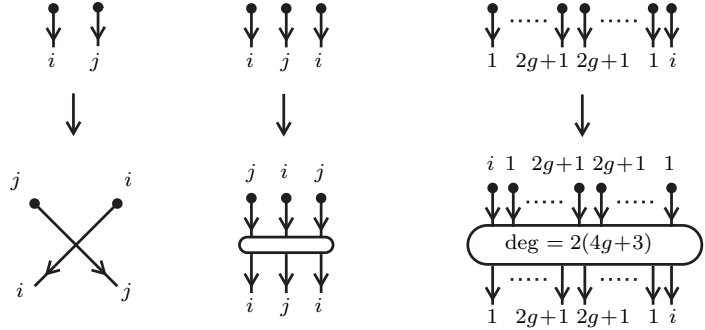}}
\put(-262,45){\footnotesize{$j$}}
\put(-219,45){\footnotesize{$i$}}
\put(-260,3){\footnotesize{$i$}}
\put(-220,3){\footnotesize{$j$}}
\put(-248,90){\footnotesize{$i$}}
\put(-232,90){\footnotesize{$j$}}
\put(-173,48){\footnotesize{$j$}}
\put(-160,48){\footnotesize{$i$}}
\put(-147,48){\footnotesize{$j$}}
\put(-173,-3){\footnotesize{$i$}}
\put(-160,-3){\footnotesize{$j$}}
\put(-147,-3){\footnotesize{$i$}}
\put(-173,90){\footnotesize{$i$}}
\put(-160,90){\footnotesize{$j$}}
\put(-147,90){\footnotesize{$i$}}
\put(-82,51){\footnotesize{$i$}}
\put(-75,51){\footnotesize{$1$}}
\put(-65,51){\footnotesize{$2g\! +\! 1$}}
\put(-42,51){\footnotesize{$2g\! +\! 1$}}
\put(-15,51){\footnotesize{$1$}}
\put(-82,-7){\footnotesize{$1$}}
\put(-71,-7){\footnotesize{$2g\! +\! 1$}}
\put(-48,-7){\footnotesize{$2g\! +\! 1$}}
\put(-21,-7){\footnotesize{$1$}}
\put(-13,-7){\footnotesize{$i$}}
\put(-75,21){\footnotesize{deg $=2(4g\! +\! 3)$}}
\put(-82,90){\footnotesize{$1$}}
\put(-71,90){\footnotesize{$2g\! +\! 1$}}
\put(-48,90){\footnotesize{$2g\! +\! 1$}}
\put(-21,90){\footnotesize{$1$}}
\put(-13,90){\footnotesize{$i$}}
\end{center}
\vspace{-0.5cm}
\caption{Chart moves}
\label{moves21}
\end{figure}

{\it Proof.} C2-move, C3-move, C4-move and their inverses 
(Figure~\ref{moves1}) 
realize H1-move, H2-move, H3-move and their inverses.  See Figure \ref{moves21}. \qed 

Now we prove Lemma~\ref{lem:change}. 

{\it Proof of Lemma~\ref{lem:change}.} 
We shall construct a chart $P_{2,h}$ by applying C2-moves, C3-moves, C4-moves and their inverses to 
$(h+1)N_0$. Such $P_{2,h}$ obviously satisfies the conditions (1) and (2). 
By Lemma~\ref{lem:hmoves}, it suffices to show that 
\begin{align*}
W'_{2,h} & = (\zeta_{2g+1},\ldots ,\zeta_1, (\zeta_{2g-2h+1},\ldots ,\zeta_{2g+1}), \ldots , (\zeta_1,\ldots ,\zeta_{2h+1}), \\
& (\zeta_1, \ldots ,\zeta_{2h})^{4h+2}, (\zeta_{2h+1},\ldots ,\zeta_1), \ldots , (\zeta_{2g+1},\ldots ,\zeta_{2g-2h+1}),\zeta_1,\ldots ,\zeta_{2g+1}) 
\end{align*}
is obtained from $W_0^{h+1}$ 
by H1-moves, H2-moves, H3-moves, their inverses, and cyclic permutations of components. 
Note that a cyclic permutation of components of a Hurwitz system does not affect a chart 
description. 
Applying H3-moves to $W_0^{h+1}=(\zeta_1, \zeta_2, \ldots, \zeta_{2g}, \zeta_{2g+1}, 
\zeta_{2g+1}, \zeta_{2g}, \ldots, \zeta_2, \zeta_1)^{2(h+1)}$, we have a Hurwitz system 
\[
((\zeta_1, \zeta_2, \ldots, \zeta_{2g}, \zeta_{2g+1})^{2(h+1)}, 
(\zeta_{2g+1}, \zeta_{2g}, \ldots, \zeta_2, \zeta_1)^{2(h+1)}). 
\]
Permuting the components of this system cyclically, we obtain 
\[
((\zeta_{2g+1}, \zeta_{2g}, \ldots, \zeta_2, \zeta_1)^{2(h+1)}, 
(\zeta_1, \zeta_2, \ldots, \zeta_{2g}, \zeta_{2g+1})^{2(h+1)}). 
\]
Applying H1-moves and H2-moves to this, we have 
\begin{align*}
& (\zeta_{2g+1},\ldots ,\zeta_1, (\zeta_{2g-2h+1},\ldots ,\zeta_{2g+1}), \ldots , (\zeta_1,\ldots ,\zeta_{2h+1}), (\zeta_{2h},\ldots , \zeta_1)^{2h+1},  \\
& (\zeta_1, \ldots ,\zeta_{2h})^{2h+1}, (\zeta_{2h+1},\ldots ,\zeta_1), \ldots , (\zeta_{2g+1},\ldots ,\zeta_{2g-2h+1}),\zeta_1,\ldots ,\zeta_{2g+1}) 
\end{align*}
by virtue of Lemma 4.6 and Lemma 4.10 of \cite{EN}. 
It follows from Lemma A.1 of \cite{Endo} that 
$(\zeta_1, \ldots, \zeta_{2h})^{2h+1}$ is obtained from 
$(\zeta_{2h}, \ldots, \zeta_1)^{2h+1}$ by H1-moves and H2-moves. 
Thus we have 
\begin{align*}
& (\zeta_{2g+1},\ldots ,\zeta_1, (\zeta_{2g-2h+1},\ldots ,\zeta_{2g+1}), \ldots , (\zeta_1,\ldots ,\zeta_{2h+1}),  \\
& (\zeta_1, \ldots ,\zeta_{2h})^{4h+2}, (\zeta_{2h+1},\ldots ,\zeta_1), \ldots , (\zeta_{2g+1},\ldots ,\zeta_{2g-2h+1}),\zeta_1,\ldots ,\zeta_{2g+1}), 
\end{align*}
which is nothing but the Hurwitz system $W'_{2,h}$. \qed


\end{document}